\documentclass[a4,12pt,multicol,makeidx]{article}
\def\origin{%p'=(p_1,...,p_m), 
  \clearpage
\vskip-\baselineskip\vskip-\topskip%
  \vbox to 0pt{\vskip-1in%
    \hbox to 0pt{\hskip-1in%
      \hbox to 0pt{\vrule width 1cm height .4pt depth 0mm\hss}%
      \vbox to 0pt{\hrule width .4pt height 0pt depth 1cm\vss}%
    \hss}%
  \vss}%%
  \vskip-\baselineskip
  \vbox to 0pt{\vskip-1in\vskip3cm%
    \hbox to 0pt{\hskip-1in\hskip3cm%
      \hbox to 0pt{\hss\vrule width 2cm height .4pt depth 0mm\hss}%
      \vbox to 0pt{\vss\hrule width .4pt height 1cm depth 1cm\vss}%
    \hss}%
  \vss}%
\vskip5mm\hskip10mm (3cm,3cm)
}%
      \def\a{\alpha}  \def\d{\delta}  \def\o{\omega}  \def\oa{\omega^1}  \def\ob{\omega^2}
\def\on{\omega^n}  \def\oN{\omega^N}  \def\l{\lambda}    \def\Gn{\Gamma^n}  \def\Ga{\Gamma^1}
\def\Gb{\Gamma^2}  \def\GN{\Gamma^N}  \def\g{\gamma}  \def\gia{\gamma_{i}^{1}}  
\def\gin{\gamma_{i}^{n}}  \def\giN{\gamma_{i}^{N}}    
\def\n{\nabla}
\def\e{\varepsilon}

\def\yy{y_i^n} \def\y{(y^1,...,y^N)}    \def\ya{y^1}  \def\yb{y^2}    \def\yN{y^N}  \def\yz{y_{\alpha}}
\def\yoa{\overline{y}_{\alpha}}  \def\yo{\overline{y}}  \def\yaa{y_{\alpha}^1}  \def\yab{y_{\alpha}^2}  \def\yan{y_{\alpha}^n}  \def\yaN{y_{\alpha}^N}
      
\def\yaan{{y}_{1,\alpha}^n}    \def\yain{{y}_{i,\alpha}^n}  \def\yadn{{y}_{d,\alpha}^n}  
  \def\yaia{{y}_{i,\alpha}^1}    \def\yao{\overline{y}^1}  \def\ybo{\overline{y}^2}  \def\yno{\overline{y}^n}
\def\yNo{\overline{y}^N}  \def\ynoa{\overline{y}^n_1}  \def\ynob{\overline{y}^n_2}    \def\ynod{\overline{y}^n_d}
      
\def\za{z_1}  \def\zb{z_2}    \def\zd{z_d}
\def\N{\nabla} \def\Nn{\nabla_{y^n}}  \def\Na{\nabla_{y^1}}  
  \def\Nin{\nabla_{y_i^n}}    
\def\R{\bf R}  \def\Rd{{\bf R}^{\rm d}}  \def\RN{{\bf R}^{\rm N}}    \def\RM{{\bf R}^{\rm M}}
\def\RdN{{\bf R}^{\rm d\times N}}    \def\Zd{{\bf Z}^{\rm d}}  \def\ZM{{\bf Z}^{\rm M}}
  \def\ZN{{\bf Z}^{\rm N}}  \def\Nat{\bf N}
  \def\T{\bf T}  \def\Td{{\bf T}^{\rm d}}  \def\TN{{\bf T}^{\rm N}}  \def\TM{{\bf T}^{\rm M}}    \def\TdN{{\bf T}^{\rm d\times N}}

  \def\ue{u_{\varepsilon}}      \def\vo{\overline{v}}
\def\vu{\underline{v}}  \def\wo{\overline{w}}  \def\wu{\underline{w}}  \def\vo{\overline{v}}  \def\vu{\underline{v}}  \def\Ho{\overline{H}}
\def\p{\partial}        \def\pd{>_{{\bf R}^d \times {\bf R}^d}}

\newcommand\sa{\sup_{\alpha\in A}}

\newenvironment{theorem}{%
\par \bigskip \it}{%
\bigskip \par}
\newenvironment{proposition}{%
\par \bigskip \it}{%
\bigskip \par}
\newenvironment{lemma}{%
\par \bigskip \it}{%
\bigskip \par}
\newenvironment{corollary}{%
\par \bigskip \it}{%
\bigskip \par}
\newenvironment{definition}{%
\par \bigskip \it}{%
\bigskip \par}

\newenvironment{remark}{%
\par \bigskip \it}{%
\bigskip \par}

\makeindex
\title{Multiscale Homogenizations for first-order Hamilton-Jacobi-Bellman equations. 
}
\author{Mariko Arisawa\\ Department of Mathematics
\\University of  Wisconsin-Madison\\
Madison, WI 53706, U.S.A.\\
}
\date{}
\begin{document}
\maketitle
\bigskip

\section{Introduction}	
 $\quad$ In this paper, we are concerned with the asymptotic behaviors of the solutions $u_{\e}(x,t)$ and $u_{\e}(x)$ of the first-order Hamilton-Jacobi-Bellman  equations (H-J-Bs in short) as $\e=(\e^1,...,\e^{n},...,\e^N)$$\in {\bf R}^{d\times N}$ tends to $0$. First, we study \\
(Evolution problem)
\begin{equation}\label{homoe}
	\frac{\p u_{\e}}{\p t}(x,t)+H(x,\frac{x}{\e^{1}},...,\frac{x}{\e^{n}},...,\frac{x}{\e^{N}},\n u_{\e})=0 \quad t\geq 0,\quad x\in \Omega={\bf R^d},
\end{equation}
$$
	u_{\e}(x,0)=u_0(x)\qquad x\in \Omega= {\bf R^d}. 
$$
(Stationary problem)
\begin{equation}\label{homos}
	\mu u_{\e}(x)+H(x,\frac{x}{\e^{1}},...,\frac{x}{\e^{n}},...,\frac{x}{\e^{N}},\n u_{\e})=0 \quad x\in \Omega,
\end{equation}
$$
	u_{\e}(x)_{|\p\Omega}=0 \qquad x\in \p\Omega,  
$$
where $\mu>0$. \\
 In (1), $\Omega={\bf R}^d$; in (2), $\Omega$ is an open bounded connected domain in ${\bf R}^d$; $\n=\n_x$;
$$
	\e=(\e^1,...,\e^n,...,\e^N)\in {\bf R}^{d\times N};
$$
$$
	\e^1=(\e^1_1,...,\e^1_1)\in {\bf R}^{d},\quad \e^n=(\e^n_1,...,\e^n_i,...,\e^n_d)\in {\bf R}^{d}\quad 2\leq \forall n\leq N;
$$
$$
	\frac{x}{\e^{n}}= (\frac{x_1}{\e^{n}_1},...,\frac{x_i}{\e^{n}_i},...,\frac{x_d}{\e^{n}_d}) \quad 1\leq \forall n\leq N;
$$
and the Hamiltonian $H$ is given by the following equation. 
\begin{equation}\label{HH}
	H(x,y^1,...,y^N,p)=\sup_{\a\in \mathcal{A}}\{
	-\left\langle b(x,y^1,...,y^N,\a),p \right\rangle_{{{\bf R^d}\times{\bf R^d}}} -g(x,y^1,...,y^N,\a)
	\},
\end{equation}
$$
	(x,y^1,...,y^N,p)\in \Omega\times {\bf R}^{d\times N}\times {\bf R^d}, \quad (y^n=(y^n_1,y^n_2,...,y^n_d)\quad 1\leq n\leq N),
$$
where $\mathcal{A}$ is a closed subset of a metric space; $b(\cdot)$ is a $d$-dimensional vector-valued function defined in $\Omega\times {\bf R}^{d\times N}\times \mathcal{A}$; $g(\cdot)$ is a real valued function defined in $\Omega\times {\bf R}^{d\times N}\times  \mathcal{A}$; and $\phi=b,g$ satisfies 
\begin{equation}\label{aasp}
	\phi(x,y^1,...,y^N,\a)\quad\hbox{is periodic in}\quad (y^1,...,y^N)\in {\bf T^{d\times N}}\quad\hbox{for}\quad \forall x\in \Omega,\forall \a\in A,
\end{equation}
\begin{equation}\label{lions}
	|\phi(x,y^1,...,y^N,\a)-\phi(x,\yao,...,\yNo,\a)|\leq L(x,\a)|(y^1,...,y^N)-(\yao,...,\yNo)|
\end{equation}
$$
	\hbox{for}\quad \forall (y^1,...,y^N), \quad \forall (\yao,...,\yNo)\in {\bf R^{d\times N}} \quad \hbox{for}\quad \forall \a\in \mathcal{A}, 
$$
where $|\cdot -\cdot|$ denotes the norm in ${\bf R^{d\times N}}$; $L(x,\a)>0$ is a constant depending on $(x,\a)\in \Omega\times \mathcal{A}$. 
The solutions $u_{\e}(x,t)$, $u_{\e}(x)$ satify (1), (2) respectively in the sense of viscosity solutions. \\

We assume that $u_0\in BUC({\bf R^d})$ and that 
\begin{equation}\label{growth}
 |H(x,y^1,...,y^N,p)|\to \infty\quad \hbox{as}\quad |p|\to \infty,\quad \hbox{uniformly for}\quad (x,y)\in \Omega\times {\bf R^{d\times N}}.
\end{equation}
Then, as in P.-L. Lions, G. Papanicolaou and S.R.S. Varadhan \cite{lpv}, L.C. Evans \cite{evans}, \cite{evansb}, $u_{\e}(x,t)$, $u_{\e}(x)$ converge uniformly on ${\bf R^{d}}\times [0,T]$ (for all $T<\infty$), on $\Omega$ respectively to some functions $\overline{u}(x,t)$, $\overline{u}(x)$ by finding the effective equtions for them: \\
(Evolutionary problem)
\begin{equation}\label{arima}
	\frac{\p \overline{u}}{\p t}(x,t)+\overline{H}(x,\n \overline{u})=0 \quad t\geq 0,\quad x\in \Omega={\bf R^d},
\end{equation}
$$
	\overline{u}(x,0)=u_0(x) \quad x\in {\bf R^d},
$$
(Stationary problem)
\begin{equation}\label{arimb}
	\mu \overline{u}(x)+\overline{H}(x,\n \overline{u})=0 \quad x\in \Omega,
\end{equation}
$$
	\overline{u}(x)_{|\p\Omega}=0 \quad x\in \p\Omega. 
$$
$\quad$ We refer the readers to A. Bensoussan, J.L. Lions and G. Papanicolaou \cite{jllions} and L. Tartar \cite{tartar} for the origin of the homogenization theory and for its various applications to physical and chemical models. A. Bensoussan, J.L. Lions and G. Papanicolaou  studied in \cite{jllions} the periodic homogenizations for linear equations. 
Later, P.-L. Lions, G. Papanicolaou, and S.R.S. Varadhan, and L.C. Evans generalized the results in \cite{jllions} to nonlinear H-J-Bs in \cite{lpv}, \cite{evans}, and \cite{evansb} by using the viscosity solutions theory. Their results correspond to  the case $N=1$ in (1) and (2) in the present paper. L. Tartar studied in \cite{tartar} the homogenization problem in a different framework than \cite{jllions} (in the framework of weakly converging functions), without periodicity assumptions. We shall study in the latter part of this paper general nonlinear homogenization problems which include the quasi periodic case and the almost periodic case. \\
By putting 
\begin{equation}\label{gin}
	\gamma^n_i=\lim_{\e\to 0}\frac{\e^{1}_{1}}{\e^{n}_{i}} \quad 1\leq i\leq d,\quad 1\leq n\leq N,
\end{equation}
we assume that the following $N$ $d\times d$ matrices $\Gamma^n$ ($1\leq n\leq N$) satisfy the Condition A below.
\begin{equation}\label{Gn}
	\Gamma^n= 
\left[ 
\begin{array}{ccc}
\gamma^n_1&\quad &O \\
\quad&\gamma^n_i &\quad \\
O&\quad &\gamma^n_d \\
\end{array} 
\right]
\end{equation}
{\bf Condition A}\\
$\gamma^n_i\neq 0,\infty$ $(1\leq \forall i\leq d, 1\leq \forall n\leq N)$, and if there exist $N$ $d\times d$ diagonal matrices with integer diagonal elements
 $Z_n$ ($1\leq n\leq N$):
$$Z_n=
\left[ 
\begin{array}{ccc}
z^n_1&\quad &O \\
\quad&z^n_i &\quad \\
O&\quad &z^n_d \\
\end{array} 
\right]\quad z^n_i \in {\bf Z}\quad 1\leq i\leq d
$$
such that all diagonal elements of $\sum_{n=1}^{N} \Gamma^n Z_n$ are integers, then $Z_n=O$ ($1\leq \forall n\leq N$).\\

Under the above assumptions, we show that the effective Hamiltonian $\overline{H}$ in (7), (8) is given by the so-called ergodic problem for H-J-Bs: 
\begin{equation}\label{ep}
	\overline{H}(x,p)=-\lim_{\l \to 0}\l w_{\l}(y^1,...,y^N)\quad \hbox{uniformly on}\quad \forall (y^1,...,y^N)\in {\bf T}^{d\times N},
\end{equation}
where for each fixed $(x,p)$, $w_{\l}$ $(\l>0)$ is the viscosity solution of the following problem. 
\begin{equation}\label{ergod}
	\l w_{\l}(y^1,...,y^N)+ H(x,y^1,...,y^N,p+\sum^{N}_{n=1}\Gamma^n \n_{y^n}w_{\l})=0
\end{equation}
$$
	\qquad\qquad\qquad(y^1,...,y^N)\in {\bf R}^{d\times N},\quad (y^n=(y^n_1,...,y^n_d)) \quad 1\leq n\leq N),
$$
\begin{equation}\label{period}
	w_{\l}\quad \hbox{is periodic in}\quad  (y^1,...,y^N)\in{\bf T}^{d\times N}, 
\end{equation}
where $\n_{y^n}=(\n_{y^n_1},...,\n_{y^n_i},...,\n_{y^n_d})$ ($1\leq n\leq N$). We derive the relationship (11) by the formal asymptotic expansion as in \cite{jllions}, \cite{lpv} and justify the formal argument by the perturbed test function method introduced in \cite{evans}. We shall see that $\overline{H}(x,p)$ defined by (11) is characterized by: for any $\delta>0$ there exist a viscosity subsolution $\underline{v}(y^1)$ and a viscosity supersolution $\overline{v}(y^1)$ of 
\begin{equation}\label{pp}
	-\overline{H}(x,p)+ H(x,\Gamma^1y^1,\Gamma^2y^1,...,\Gamma^Ny^1,p+\n_{y^1}\underline{v}(y^1))\leq \delta,\quad y^1\in {\bf R}^d,
\end{equation}
\begin{equation}\label{ll}
	-\overline{H}(x,p)+ H(x,\Gamma^1y^1,\Gamma^2y^1,...,\Gamma^Ny^1,p+\n_{y^1}\overline{v}(y^1))\geq -\delta,\quad y^1\in {\bf R}^d. 
\end{equation}
The functions $\underline{v}(y^1)$ and $\overline{v}(y^1)$ serve as "approximated" correctors in the perturbed test function method. We shall call (14)-(15) the approximated cell problem. We remark that although $\overline{H}(x,p)$ is defined by the ergodic problem (12)-(13) in ${\bf T}^{d\times N}$, it is characterized by the approximated cell problem (14)-(15) in ${\bf R^d}$.\\

 The existence of the limit on the right-hand side of (11) is not trivial. Is is closely related to the ergodicity of the controlled deterministic system $\overline{Y}^{\a}(t)$
 on the torus ${\bf T}^{d\times N}$, which is the embedding of the following system $Y^{\a}(t)$$=(y_{\a}^1(t),y_{\a}^2(t),...,y_{\a}^N(t))$ in ${\bf R}^{d\times N}$: 
\begin{equation}\label{system}
	\frac{d}{dt}(y^1_{\a}(t),...,y^N_{\a}(t))=\overline{B}(Y^{\a}(t),\a(t)) \quad t\geq 0,
\end{equation}
$$
	(y^1_{\a}(0),...,y^N_{\a}(0))=(y^1,...,y^N)\qquad (y^1,...,y^N)\in {\bf R}^{d\times N}, 
$$
where $\a(t)$ is a measurable function from $[0,\infty)$ to $\mathcal{A}$; and 
$$
	\overline{B}(y,\a)=(\Gamma^1b(x,y^1,...,y^N,\a),\Gamma^2b(x,y^1,...,y^N,\a),...,\Gamma^Nb(x,y^1,...,y^N,\a)),
$$
for $y=(y^1,...,y^N)$ (remark that $\Gamma^1=I$). To see the relationship between $\overline{Y}^{\a}(t)$ and (12)-(13), we shall rewrite (12) to 

\begin{equation}\label{rew}
	\l w_{\l}+ \sup_{\a\in \mathcal{A}}\{-\left\langle \overline{B}(y,\a),\n_y w_{\l}  \right\rangle_{{{\bf R^d}\times{\bf R^d}}}
	-\left\langle b(x,y,\a),p \right\rangle_{{{\bf R^d}\times{\bf R^d}}} 
\end{equation}
$$
	-g(x,y^1,...,y^N,\a)
	\}=0 \quad y=(y^1,...,y^N)\in {\bf R}^{d\times N}. 
$$
We find 
$$
	\lim_{\l\to 0} \l w_{\l}(y^1,...,y^N)= 
$$
\begin{equation}\label{long}
	\lim_{\l\to 0}\inf_{\a(\cdot)\in \mathcal{A}}\l\int_{0}^{\infty} e^{-\l t} \{
	g(x,\overline{Y}_{\a}(t),\a(t))+ \left\langle b(x,\overline{Y}_{\a}(t) ,\a(t)),p \right\rangle_{{{\bf R^d}\times{\bf R^d}}} 
	\}dt
\end{equation}
$$
\qquad\qquad\qquad\qquad\qquad\qquad\qquad\qquad\qquad\forall (y^1,...,y^N)\in {\bf T}^{d\times N}. 
$$
The left-hand side's limit in (18) is the long time average of the potential over the system $\overline{Y}^{\a}(t)$. We refer the reader for the ergodic problems for H-J-Bs to J.M. Lasry \cite{lasry}, P.-L. Lions \cite{lpn}, M. Arisawa \cite{infinite}, \cite{arisawaa}, \cite{arisawab}, M. Arisawa and P.-L. Lions \cite{apl}, and M. Robin \cite{robin}.\\

It is worth remarking here that the Condition A in (16) plays a similar role to the necessary and sufficient condition for the ergodic parallel transformations in tori ${\bf T}^N$'s: 
\begin{equation}\label{jacobi}
	T_i\quad :\quad x\in {\bf T^N}\to x+(\gamma^1_i,\gamma^2_i,...,\gamma^N_i)\in {\bf T^N} \quad 1\leq i\leq d.
\end{equation}
See for example V. Arnold and A. Avez \cite{arnold}, in which the Condition A appears in the Jacobi's Lemma. Recently, M. Arisawa and P.-L. Lions studied in \cite{apl} 
a general necessary and sufficient condition - called non-resonance condition - for the ergodic problem for H-J-Bs with constant coefficients. The condition A is nothing less than the non-resonance condition in \cite{apl} in the case of constant coefficients.\\

Next, we study some more general homogenization problems for the following H-J-Bs. \\

(Evolution problem)
\begin{equation}\label{ahomoe}
	\frac{\p {u}_{\e}}{\p t}(x,t)+{F}(x,\frac{x}{\e},\n {u}_{\e})=0 \quad t\geq 0,\quad x\in \Omega={\bf R^d},.
\end{equation}
$$
	{u}_{\e}(x,0)=u_0(x) \quad x\in \Omega={\bf R^d},
$$
(Stationary problem)
\begin{equation}\label{ahomos}
	\mu {u}_{\e}(x)+  {F}(x,\frac{x}{\e},\n {u}_{\e})  =0 \quad x\in \Omega,.
\end{equation}
$$
	{u}_{\e}(x)_{|\p\Omega}=0 \quad x\in \p\Omega,
$$
where $\mu>0$. \\
Here, 
$$
	\e=(\e_1,\e_2,...,\e_d)\in {\bf R}^d,\quad \frac{x}{\e}=(\frac{x_1}{\e_1},\frac{x_2}{\e_2},...,\frac{x_d}{\e_d})\in {\bf R^d}. 
$$
We consider the following three classes of Hamiltonians.\\

\leftline{\bf Definition 1.1}
\begin{definition}
We say that the Hamiltonian $F(x,y,p)$ defined in $\Omega\times {\bf R}^{d}\times {\bf R}^{d}$ is quasi periodic in $y\in {\bf R^d}$ if it satisfies (6) and 
\begin{equation}\label{gap}
	F(x,y,p)=\sum_{n=1}^{N}F^{n}(x,y,p) \quad (x,y,p)\in \Omega\times {\bf R}^{d}\times {\bf R}^{d}, 
\end{equation}
where 
$$
	F^{n}(x,y,p)=\sup_{\a\in \mathcal{A}}\{
	-\left\langle b^n(x,y,\a),p \right\rangle_{{{\bf R^d}\times{\bf R^d}}} -g^n(x,y,\a)
	\} \quad (x,y,p)\in \Omega\times{\bf R^d}\times{\bf R^d},
$$
where $b^n$ is a $d$-dimensional vector valued function; $g^n$ is a real valued function; for each $(x,p)$, $\phi^n(y)=b^n$, $g^n(x,y,p)$ ($1\leq n\leq N$) satisfy 
\begin{equation}\label{quasi}
	\phi^n(y) \quad \hbox{is periodic in}\quad y\in \Pi_{1\leq i\leq d}[0,T^n_{i}], 
\end{equation}
and 
\begin{equation}\label{la}
	|\phi^n(x,y,\a)- \phi^n(x,y',\a)|\leq L(x,\a)|y-y'| \quad \forall y,\quad y'\in {\bf R^d},
\end{equation}
where $L(x,\a)$ depends only on $(x,\a)$. 
\end{definition}

\leftline{\bf Definition 1.2}
\begin{definition}
We say that the Hamiltonian $F(x,y,p)$ defined in $\Omega\times {\bf R^d}\times {\bf R^d}$ is in the class $\mathcal{B}_0$ Hamiltonian if 
\begin{equation}\label{unip}
	F(x,y,p)=\lim_{N\to \infty} F^N(x,y,p) \quad \hbox{uniformly in}\quad (x,y,p)\in \Omega\times{\bf R^d}\times{\bf R^d},
\end{equation}
where $F^N(x,y,p)$ ($N\in {\bf N}$) is uniformly bounded, quasi periodic in $y\in {\bf R^d}$ with constants $L(x,\a)$ in (24) independent on $N\in {\bf N}$, and satisfies (6) uniformly in $N\in {\bf N}$. 
\end{definition}

\leftline{\bf Definition 1.3}
\begin{definition}
  We say that the Hamiltonian $F(x,y,p)$ defined in $\Omega\times {\bf R^d}\times {\bf R^d}$ is in the class $\mathcal{B}_1$ Hamiltonian if it satisfies (6) and 
\begin{equation}\label{gaa}
	F(x,y,p)=\sup_{\a\in \mathcal{A}}\{
	-\left\langle \a,p \right\rangle_{{{\bf R^d}\times{\bf R^d}}} -V(x,y,\a)
	\} \quad (x,y,p)\in \Omega\times{\bf R^d}\times{\bf R^d},
\end{equation}
where for each $(x,\a)$, $V$ is bounded and satisfies 
\begin{equation}\label{general}
	\lim_{T\to \infty} \frac{1}{T} \int_{0}^{T} 
	V(x,y+tz,\a) dt=c(x,\a)\quad \hbox{uniformly in}\quad y, \quad z  (|z|=1) \in {\bf R}^d, 
\end{equation}
where for each $(x,\a)$, there exists a constant $c(x,p)$ such that 
\begin{equation}\label{cancel}
	c(x,\a)+ \left\langle x,\a \right\rangle_{{\bf R^d}\times{\bf R^d}}-c(x,p)\leq 0 \quad \forall \a\in \mathcal{A},\quad |\a|=1.
\end{equation}
\end{definition}

It is well known that if a function $V(y)$ is almost periodic in $y\in {\bf R^d}$ in the sense of H. Bohr (\cite{bohr}, see Definition 1.4 below), then $V$ satisfies (27). 
We shall give some examples of Hamiltonians in the class $\mathcal{B}_1$ below in this introduction. (The examples of quasi periodic and the class $\mathcal{B}_0$ Hamiltonians are obvious from the definitions.)\\

\leftline{\bf Definition 1.4}
\begin{definition}
Let $f(y)$ be a function defined in ${\bf R^d}$. The function $f(y)$ is almost periodic in the sense of H. Bohr if and only if the set of functions 
\begin{equation}\label{Bohr}
	\{ f(y+tz)\quad|\quad \forall t\in {\bf R},\quad \forall z\in {\bf R^d}  \}
\end{equation}
is relatively compact in the space of bounded continuous functions on ${\bf R^d}$ with the uniform norm $||f||$$=\sup_{y\in {\bf R^d}}|f(y)|$. 
\end{definition}

We assume that $F$ is either quasi periodic, or in the class $\mathcal{B}_0$, or in the class $\mathcal{B}_1$. We also assume that $u_0$ is $BUC({\bf R^d})$. 
 Under these assumptions, we show that $u_{\e}(x,t)$, $u_{\e}(x)$ of (20), (21) converge to $\overline{u}(x,t)$, $\overline{u}(x)$ uniformly in $\Omega\times [0,T]$, 
$\Omega$ respectively the solutions of the following problems.\\
(Evolutionary problem)
\begin{equation}\label{aarima}
	\frac{\p \overline{u}}{\p t}(x,t)+\overline{F}(x,\n \overline{u})=0 \quad t\geq 0,\quad x\in \Omega={\bf R^d},
\end{equation}
$$
	\overline{u}(x,0)=u_0(x) \quad x\in \Omega={\bf R^d},
$$
(Stationary problem)
\begin{equation}\label{aarimb}
	\mu \overline{u}(x)+\overline{F}(x,\n \overline{u})=0 \quad x\in \Omega={\bf R^d},
\end{equation}
$$
	\overline{u}(x)_{|\p\Omega}=0 \quad x\in \p\Omega. 
$$

For the quasi periodic Hamiltonian $F$, the effective Hamiltonian $\overline{F}$ is given by the ergodic problem for H-J-Bs: 
\begin{equation}\label{aep}
	\overline{F}(x,p)=-\lim_{\l\to 0} \l v_{\l}(y)\quad \hbox{uniformly in}\quad y\in {\bf R^d}, 
\end{equation}
where for each $(x,p)$, $v_{\l}$ ($\l>0$) is the viscosity solution of 
\begin{equation}\label{aergod}
	\l v_{\l}(y)+F(x,y,p+\n_y v_{\l}(y))=0 \quad y\in {\bf R^d},
\end{equation}
which grows at most linearly at infinity.\\

The effective Hamiltonian $\overline{F}(x,p)$ is characterized by: for any $\d>0$, there exist a bounded uniformly continuous viscosity subsolution $\underline{v}$ and 
a bounded uniformly continuous viscosity supersolution $\overline{v}$ of 
$$
	-\overline{F}(x,p)+ F(x,y,p+\n_y \underline{v}(y)) \leq \d\quad y\in {\bf R^d}, 
$$
\begin{equation}\label{mm}
	\quad
\end{equation}
$$
	-\overline{F}(x,p)+ F(x,y,p+\n_y \overline{v}(y)) \geq -\d\quad y\in {\bf R^d}. 
$$
For the class $\mathcal{B}_0$ Hamiltonian $F$, the effective Hamiltonian $\overline{F}$ is given by: 
\begin{equation}\label{alep}
	\overline{F}(x,p)=-\lim_{N \to \infty} \overline{F}^N (x,p), 
\end{equation}
where $\{\overline{F}^N(x,p)\}_{N\in {\bf N}}$ are defined by (32) for $F^N$ the quasi periodic Hamiltonians which approximate $F$.\\
For the class $\mathcal{B}_1$ Hamiltonian $F$, the effective Hamiltonian $\overline{F}$ is given by: 
\begin{equation}\label{glep}
	\overline{F}(x,p)=-\lim_{\l'\to 0} \l' v_{\l'}(y)\quad \hbox{locally uniformly in}\quad y\in {\bf R^d}, 
\end{equation}
where for each $(x,p)$, $v_{\l}$ ($\l>0$) is the viscosity solution of (33) which grows at most linearly, and $\{v_{\l'}\}$ is a subsequence of $\{v_{\l}\}$. 
The effective Hamiltonian $\overline{F}(x,p)$ is characterized  by: there exists a uniformly continuous viscosity solution $v$ of (34) such that 
\begin{equation}\label{co}
	\lim_{|y|\to \infty} \frac{v(y)}{|y|}=0.
\end{equation}

 As we shall see later, our quasi periodic homogenizations are special case of (1), (2), and $\overline{H}(x,p)$$=\overline{F}(x,p)$ in such cases. We remark that the same effective Hamiltonian $\overline{H}(x,p)$$=\overline{F}(x,p)$ is derived by two different ergodic problems: (12)-(13) in ${\bf T^{d\times N}}$ and (33) in ${\bf R^d}$. \\

Here, we shall give some examples which illustrate the results in this paper.\\

\leftline{\bf Example 1.1}
\begin{definition}
 Let 
\begin{equation}\label{qual}
	H(x,y^1,...,y^N,p)=a(x,y^1,...,y^N)|p|^2-V(x,y^1,...,y^N),
\end{equation}
for $(x,y^1,...,y^N,p)\in \Omega\times {\bf R^{d\times N}}\times {\bf R^d}$, where $a(\cdot)\geq a>0$, $V(\cdot)$ satisfy (4), (5). Consider (1), (2) with (38). 
(Remark that we can rewrite (38) to (3)). Then, the solution $u_{\e}(x,t)$, $u_{\e}(x)$ converge uniformly to the solutions $\overline{u}(x,t)$, $\overline{u}(x)$ of 
(7), (8), where the effective Hamiltonian $\overline{H}$ is given by (11). 
\end{definition}

\leftline{\bf Example 1.2}
\begin{definition}
 Let 
\begin{equation}\label{finite}
	H(x,y^1,...,y^N,p)=a(x,y^1)|p|-V(x,y^1,...,y^N),
\end{equation}
for $(x,y^1,...,y^N,p)\in \Omega\times {\bf R^{d\times N}}\times {\bf R^d}$, where $a(\cdot)$ is a function of $(x,y^1)\in \Omega\times {\bf R^d}$ such that 
$a(\cdot)\geq a>0$, and $V(\cdot)$ satisfies  (4), (5). Consider (1), (2) with (39). 
(Remark that we can rewrite (39) to (3)). Then, the solution $u_{\e}(x,t)$, $u_{\e}(x)$ converge uniformly to the solutions $\overline{u}(x,t)$, $\overline{u}(x)$ of 
(7), (8), where the effective Hamiltonian $\overline{H}$ is given by (11). 
\end{definition}

\leftline{\bf Example 1.3}
\begin{definition}
 Let 
\begin{equation}\label{dyna}
	F(x,y,p)=|p|^2-V(x,y),
\end{equation}
where for each $x\in \Omega$, $V(x,y)$ is almost periodic in the sense of H. Bohr. It is easy to see that $F$ is a class $\mathcal{B}_1$ Hamiltonian. Consider 
(20), (21) with (40). Then, the solution $u_{\e}(x,t)$, $u_{\e}(x)$ converge uniformly to the solutions $\overline{u}(x,t)$, $\overline{u}(x)$ of 
(30), (31), where the effective Hamiltonian $\overline{F}$ is given by (36). 
\end{definition}

\leftline{\bf Example 1.4}
\begin{definition}
 Let $d=1$, 
\begin{equation}\label{almost}
	F(x,y,p)=|p|-V(y) \quad (x,y,p)\in \Omega\times {\bf R}\times {\bf R},
\end{equation}
$$
	V(y)=\quad |y|\quad |y|\leq 1;\quad =1\quad |y|\geq 1.
$$
It is easy to see that $V(\cdot)$ satisfy (27) and that $F$ is a class $\mathcal{B}_1$ function. Consider (20), (21) with (41). Then, the solution $u_{\e}(x,t)$, $u_{\e}(x)$ converge uniformly to the solutions $\overline{u}(x,t)$, $\overline{u}(x)$ of 
(30), (31), where the effective Hamiltonian $\overline{F}$ is given by (36). 
\end{definition}

We shall give the proofs of Examples 1.1-1.4 below in this paper.\\

Now, we shall give the plan of this paper. In \S 2, we study the ergodic problem for H-J-Bs: (12)-(13) in ${\bf T^{d\times N}}$. We shall prove the existence of (11), by studying the controllability of the system (16). In \S 3, the approximated cell problem (14)-(15) is solved by using the result in \S 2. This section describes a important relationship between the ergodic problems in ${\bf T^{d\times N}}$ and 
in ${\bf T^{d}}$ for quasi periodic homogenizations. In \S 4, we prove the effective equations (7), (8) rigorously. This is the first main result. In \S 5, we study the general homogenizations. This is the second main results. In \S 6, we describe the derivation of (11) by the formal asymptotic expansions argument.\\

Before giving notational remarks and concluding this introduction, let us give some remarks concerning with this paper. Although we assume that $\gamma^n_i\neq 0$, $\infty$, we can treat the cases when $\gamma^n_i=0$, $\infty$ for some $1\leq i\leq d$, $1\leq n\leq N$. For such cases, we obtain the effective Hamiltonian 
by using iterations of homogenizations in addition to the results in this paper. For the iteration of homogenizations, see \cite{avellaneda},   \cite{fletcher},  \cite{jllions},  etc. As is easily seen, we can relax some minor assumptions. For example, while we assumed in this paper that $\e^1\in {\bf R^d}$ (the first $d$-dimensional entry in $\e\in {\bf R^{d\times N}}$) was $\e^1=(\e^1_1,...,\e^1_1)$, we can consider the homogenization with $\e^1$ with $d$ different entries in a similar way. \\
We shall study the multiscale homogenization for second-order H-J-Bs in the forthcoming paper.\\ 

Throughout this paper, we use the notion of viscosity solutions which was introduced by M.G. Crandall-P.-L. Lions in \cite{crandallpl}. We refer the reader M.G. Crandall-H. Ishii-P.-L. Lions \cite{crandall}, W.H. Fleming-H.M. Soner \cite{fleming}.\\

We denote by ${\bf R,Q,Z,N}$ the sets of real, rational, integer, natural numbers, by ${\bf R^M,T^M}$ the $M$-dimensional Euclidian space and $M$-dimensional 
torus respectively. The distance between two points $y,y'$$\in {\bf R^M}$ is denoted by $|y-y'|$ for any dimensional number $M$. On 
the other hand, the distance between
 two points  $y,y'\in$ $\T^M$ is denoted by $d_M (y,y')$. In $\R^M$, for $x\in \R^M$
$U_r(x)$ denotes the ball of radius $r>0$, centered at $x$.
We denote 
$$
	{\bf Z}^M=\{(z_1,z_2,...,z_M)|\quad z_i\in {\bf Z},\quad 1\leq \forall i\leq M \},
$$
$$
	{\bf N}^M=\{(z_1,z_2,...,z_M)|\quad z_i\in {\bf N},\quad 1\leq \forall i\leq M \}.
$$
We  use right-bottom indices $x_i$ $(1\leq i\leq d)$ to represent the $i$-th entry in  $d$- dimensional spaces; we use right-upper indices   $y^n$ $(1\leq n\leq N)$ to represent the $n$-th rescaled variable; and we use the both indices  $y^n_i$ to represent the $i$-th entry of the $n$-th rescaled variable. The rule of the indices notations is invariant throughout this
paper. (When we need to write a power of a number, a specific remark will be added in each occation.)\\

\bigskip

	The author began to be interested in the homogenization problems
from the works and lectures of Professors P.-L. Lions, L.C. Evans, P.E. Souganidis, and 
L. Tartar. The discussions with them were  very inspiring and helpful. 
She would like to sincerely thank to them.

\bigskip

\section{Ergodic problem in $\TdN$}

	We study the problem ($\ref{ergod}$)-($\ref{period}$), and show the ergodic convergence property
\begin{equation}\label{wl}
	\lim_{\l \to 0}\l w_{\l}\y \to \exists d_{x,p}^{!}\qquad 
	\hbox{uniformly in}\quad \y\in \TdN,
\end{equation}
where $d_{x,p}$ is a constant. We give some sufficient conditions for ($\ref{wl}$) in terms of controllabilities of the system $\overline{Y}_{\alpha}(t)$ in $\TdN$. First, we shall recall the notions of controllabilities introduced in $\cite{arisawab}$.

\bigskip

\leftline{\bf Definition 2.1}
\begin{definition}
	Let $c(y,\a)$ be a $\RM$-valued function defined in $(y,\a)\in \RM \times A$, and
 consider the controlled system  defined by the following ordinary
differential equation for each control $\a(\cdot)$:
\begin{equation}\label{controllable}
	\frac{d}{dt}\yz (t)=c(\yz (t),\a (t))  \qquad t\geq 0,\qquad
	\yz (0)=y.
\end{equation}
$$\quad$$ 
\begin{enumerate}
\item A point $y\in \RM$ is said to be exactly controllable to a point $y'\in \RM$
in the system ($\ref{controllable}$)  if there exist a control $\a (\cdot)$
and $T(y,y')>0$ such that $\yz(T(y,y'))=y'$.

\item	A point $y\in \RM$ is said to be approximately controllable to a point 
$y'\in \Rd$ in the system ($\ref{controllable}$) in $\RM$
with the estimate $\Delta (\d; y,y')$,
 if for any $\d>0$ there exist a control $\a (\cdot)$
and $T(\d;y,y')>0$ such that 
$|y_{\alpha} (T(\d ;y,y'))-y'|<\d $, 
$T(\d;y,y')<\Delta (\d; y,y')$.
\end{enumerate}
\end{definition}

\leftline{\bf Definition 2.2}
\begin{definition}
	Let $c(y,\a)$  in Definition 2.1 be  periodic in $y\in \TM$.
Let $\yoa (t)$ be the embedding of the system $y_{\a}(t)$ in ($\ref{controllable}$) 
to $\TM$:
\begin{equation}\label{embed}
	\yoa (t) \equiv y_{\a}(t) \quad \pmod{1^M} \quad \in \TM\quad \forall t\geq 0,
	\qquad \yoa (0)=\yo \equiv y \quad \pmod{1^M}.
\end{equation}

\begin{enumerate}
\item	A point $\yo \in \TM$ is said to be exactly controllable to a point $\yo'\in \TM$
in the system ($\ref{embed}$) if there exist a control $\a (\cdot)$
and $T(\yo,\yo')>0$ such that $\yoa(T(\yo,\yo'))=\yo'$.

\item	A point $\yo\in \TM$ is said to be approximately controllable to a point 
$\yo'\in \TM$ in the system ($\ref{embed}$)
with the estimate $\Delta (\d; \yo,\yo')$,
 if for any $\d>0$ there exist a control $\a (\cdot)$
and $T(\d;\yo,\yo')>0$ such that $d_M(\yoa(T(\d;\yo,\yo')),\yo')<\d$, 
$T(\d;\yo,\yo')<\Delta(\d; \yo,\yo')$.
\end{enumerate}
\end{definition}

\leftline{\bf Remark 2.1}
\begin{remark}
	Let $\yoa (t)$ be the embedding of $y_{\a}(t)$  in $\RM$ to $\TM$, and assume that 
$y=\yo$ $\in \TM$. We can study the controllability of the embedded system in $\TM$ by the controllability of the system in $\RM$.
\begin{enumerate}
\item	A  point  $\yo \in \TM$ is  exactly controllable to a point $\yo'\in \TM$
in the system ($\ref{embed}$) if and only if there exist $z\in \ZM$, a control $\a(\cdot)$
and $T(\yo,\yo')>0$ such that $\yz(T(\yo,\yo'))=\yo'+z$ in $\RM$.

\item	A point $\yo\in \TM$ is approximately controllable to a point $\yo'\in \TM$ in the system ($\ref{embed}$) with the estimate $\Delta(\d; \yo,\yo')$
 if and only if for any $\d>0$ there exist $z_{\d}\in \Zd$, a control $\a (\cdot)$
and $T(\d;\yo,\yo')>0$ such that 
$|\yz(T(\d;\yo,\yo'))-(\yo'+z_{\d})|<\d$,  $T(\d;\yo,\yo')<\Delta (\d; \yo,\yo')$.
\end{enumerate}
\end{remark}

In  \cite{arisawab}, we gave some sufficient conditions for the ergodic convergence property by using the notions of controllabilities.
Here, we   recall a part of the results in  \cite{arisawab} in
a slightly modified form. We consider
\begin{equation}\label{phin}
	\l u_{\l}(y) + \sa
	 \{-<c(y,\a),\N u_{\l}(y) \pd -h(y,\a)\}=0 \qquad y\in \TM,
\end{equation}
$$
	u_{\l}(y) \quad \hbox{is periodic in}  \quad y\in \TM,
$$
where $c(\cdot)$, $h(\cdot)$ are periodic, and Lipschitz continuous in $y$ for each $\a \in A$. 
We say that the system ($\ref{controllable}$) is  Lipschitz continuous if

\begin{equation}\label{Lip}
	|y_{\a}(t)-y'_{\a}(t)|\leq C|y-y'|\qquad \forall y,\quad y'\in \RM, \quad 		\forall \a \in A, \quad \forall t\geq 0,
\end{equation}
where $C>0$ is a constant.
We say that system ($\ref{controllable}$) is  partially Lipschitz continuous if

\begin{equation}\label{PLip}
	|y_{\a}(t)-y'_{\a}(t)|\leq C(y^1,...,y^m)|y-y'|\qquad \forall \a \in A, \quad \forall t\geq 0,
\end{equation}
$$
	\forall y=(y^1,y^2,...,y^m,y^{m+1}...,y^M),\quad 
	y'=(y^1,y^2,...,y^m,{y'}^{m+1},...,{y'}^{M})\in \RM, 
$$
where $C(y^1,...,y^m)>0$  depends only on $(y^1,...,y^m)$. We introduce the following conditions
for $h(\cdot)$
\begin{equation}\label{bdd}
	|h(y,\a)|\leq C_{\alpha}\qquad \forall y\in \RM,
\end{equation}
where $C_{\a}$ is a constant depending on $\a$, and
\begin{equation}\label{power} 
	|h(y,\a)-h(y',\a)|\leq 
	C_0(\inf_{y\in \Omega}h(y,\a)+C_1)|y-y'| \qquad  \forall y,\quad y'\in \RM, \quad 		\forall \a \in A,
\end{equation}
where $C_0$, $C_1>0$ are constants.

\bigskip

\leftline{\bf Theorem A. (\cite{arisawab})}
\begin{theorem}
Assume that either one of the following three cases holds. \\

(i) (Uniform approximated controllability) 
The functions $c(\cdot)$, $h(\cdot)$ in ($\ref{phin}$) are uniformly bounded in
 $\TM \times A$, and Lipschtz continuous  in $\TM$ uniformly in $\alpha \in A$.
There exist  $\gamma\in [0,1)$, $C>0$ such that 
any $\yo\in \TM$ is approximately controllable to any $\yo' \in \TM$ with the estimate
$\Delta(\delta;\yo,\yo')$ such that
\begin{equation}\label{estimate}
	\Delta(\delta;\yo,\yo')\leq C(-\log \d)^{\gamma}, \qquad \forall \delta>0,\quad
	\forall \yo,\quad \yo'\in \TM.
\end{equation}\\

(ii) (Approximated controllability in the Lipshitz continuous system)
Let the system be Lipschitz continuous (($\ref{Lip}$)). Any $\yo\in \TM$ is approximately controllable to any $\yo' \in \TM$ with the estimate $\Delta(\delta;\yo,\yo')$, where $\Delta$ may diverge as $\d$ goes to $0$.\\

(iii) (Approximated controllability in the partially Lipschitz continuous system)
Let the system be partially Lipschitz continuous (($\ref{PLip}$)). 
There exists a finite number $T>0$ such that  for any ${y'}^{1}$, ${y'}^{2}$,..., ${y'}^{m}$ $\in [0,1]$, and for any $\yo \in\TM$, there exists a point  $\yo'\in \TM$ such that 
$$
	\yo'=({y'}^{1},...,{y'}^{m}, {y'}^{m+1},...,{y'}^{M}),
$$
and $\yo$ 
is exactly controllable to $\yo'$  with $T(\yo,\yo')<T$. Moreover, if $\yo$, 
$\yo'$ $\in \TM$ have the same first $m$ entries, then $\yo$ is approximately controllable to $\yo'$  with the estimate $\Delta(\delta;\yo,\yo')$, where $\Delta$ may diverge as $\d$ goes to $0$.\\

	Consider the problem ($\ref{phin}$). Then, there exists a unique constant $r\in {\bf R}$ such that
$$
	\lim_{\l \to 0}\l u_{\l} (y)= r  \qquad \hbox{uniformly in}\quad y\in \TM.
$$
\end{theorem}

\par

We shall give a brief proof of Theorem A in the end of this section. 
In view of (\ref{rew}), we shall apply Theorem A to study the ergodic convergence (\ref{wl}) by examining the 
controllability of the system $\overline{Y}_{\a}(t)$ in $\TdN$ defined by ($\ref{system}$). From the definition of $\overline{B}$ in ($\ref{system}$), if we write

$$
	Y_{\a}(t)=(y^1_{\a}(t),...,y^n_{\a}(t),...,y^N_{\a}(t)),
$$
where
$$
	\yan(t)=(\yaan(t),...,\yain(t),...,\yadn(t)) 
	\qquad 1\leq i\leq d,\quad 1\leq n\leq N,
$$
we have

\begin{equation}\label{zoo}
	\yain(t)= \yy+ \gin (\yaia(t)-y_i^1)\qquad t\geq 0,
	\quad 1\leq i\leq d,\quad 1\leq n\leq N.
\end{equation}
$$\quad$$
Therefore, we see that 
the controllability of ($\ref{system}$) is determined   by the 
following  O.D.E. with respect to the variable $y^1\in \Rd$:
\begin{equation}\label{first}
	\frac{d}{dt}\yaa(t)=
	 b(x, \yaa(t),\Gb(\yaa(t)-\ya)+\yb,...,\GN(\yaa(t)-\ya)+y^N,\a(t))
	\quad t\geq 0,
\end{equation}
$$
	\yaa(0)=\ya.
$$

The following result describes the relationship between the  controllabity of 
$\yaa(t)$ in $\Rd$ and that of  $\overline{Y}_{\a}(t)$ in $\TdN$.\\

\bigskip

\leftline{\bf Proposition 2.1}
\begin{proposition}
Assume that for any $\yb,...,\yN$ $\in \Td$, any $y^1$ $\in \Rd$ is approximately controllable to 
any ${y'}^1$ $\in \Rd$ in the system ($\ref{first}$) with the estimate
$\Delta_0(\delta;y^1,{y'}^{1})$. Then, any $\yo\in \TdN$ is 
approximately controllable to any $\yo'\in \TdN$ 
 in the system $\overline{Y}_{\a}(t)$ in $\TdN$ defined by ($\ref{system}$) with the 
estimate
$\Delta(\delta;\yo,\yo')$. Moreover, if for each $\d>0$, for $\forall \rho>0$ there exists  $\Delta_0$ such that $\Delta_0<\rho$, then  for $\forall \rho>0$ there exists $\Delta$ such that $\Delta<\rho$.
\end{proposition}

\bigskip

We may now state the main result in this section.

\bigskip

\leftline{\bf Theorem 2.2}
\begin{theorem}
Let 
$$
	h_{x,p}(y,\a)=<b(x,y,\a),p\pd + g(x,y,\a) \qquad (y,\a)\in 	\TdN\times A, 
$$
and consider ($\ref{ergod}$)-($\ref{period}$).
If either one of the following three cases holds, then for any $(x,p)\in \Omega\times \Rd$,  ($\ref{wl}$) holds.\\

(i) The function $h_{x,p}$ is bounded and Lipschitz continuous in $y\in \TdN$ uniformly in 
$\a\in A$. For any $\yb,...,\yN$ $\in \Td$, any $y^1$ $\in \Rd$ is approximately controllable to any ${y'}^{1}$ $\in \Rd$  in the system ($\ref{first}$) with the estimate
$\Delta_0(\delta;y^1,{y'}^{1})$, where for each $\d>0$, for $\forall \rho>0$ there exists 
$\Delta_0<\rho$.\\

(ii) The function $h_{x,p}$ satisfies ($\ref{bdd}$), ($\ref{power}$). 
 The system ($\ref{first}$) is Lipschitz continuous uniformly in $\yb,...,\yN$ $\in \Rd$.    For any $\yb,...,\yN$ $\in \Td$, any $y^1$ $\in \Rd$ is approximately controllable to 
any ${y'}^{1}$ $\in \Rd$ in the system ($\ref{first}$) with the estimate $\Delta_0(\delta;y^1,{y'}^{1})$, which may diverge as $\d$ goes to $0$.\\

(iii) The function $h_{x,p}$ satisfies ($\ref{bdd}$), ($\ref{power}$). 
Let  $y^1_{\a}(t)$ be the solution of 
$$
	\frac{d}{dt}\yaa(t)=
	 b(x, \yaa(t),\Gb(\yaa(t)-\ya)+\yb,...,\GN(\yaa(t)-\ya)+y^N,\a(t))
	\quad t\geq 0,
$$
$$
	\yaa(0)=\ya.
$$
Let  ${y'}^1_{\a}(t)$  be the solution of  
$$
	\frac{d}{dt}{y'}^1_{\a}(t)=
 	b(x,{y'}^1_{\a}(t),
	\Gb({y'}^1_{\a}(t)-\ya)+{y'}^2,...,\GN({y'}^1_{\a}(t)-\ya)+{y'}^N,\a(t))
	\quad t\geq 0,
$$
$$
	\yaa(0)=\ya.
$$
Then,  the following holds
$$
	|y^1_{\a}(t)-{y'}^1_{\a}(t)|\leq |(\ya,\yb,...,\yN)-({y}^1,{y'}^2,...,{y'}^N)|.
$$ 
For any $\yb,...,\yN$ $\in \Rd$, any $y^1$ $\in \Rd$ is exactly controllable to 
any ${y'}^{1}$ $\in \Rd$ in the system ($\ref{first}$).

\end{theorem}

\bigskip

\leftline{\bf Remark 2.2}
\begin{remark}
	The constant $d_{x,p}$ is characterized by: for any $\d>0$ there exist a viscosity subsolution $\wu$ and a viscosity supersolution $\wo$ of

\begin{equation}\label{ari}
	d_{x,p} + H(x,y^1,...,y^N,p+\Sigma_{n=1}^N \Gn \Nn \wu) \leq \d,
\end{equation}
$$\quad$$
\begin{equation}\label{sawa}
	d_{x,p} + H(x,y^1,...,y^N,p+\Sigma_{n=1}^N \Gn \Nn \wo) \geq -\d,
\end{equation}
$$\quad$$
\begin{equation}\label{ppp}
	\wo, \quad \wu \quad \hbox{are periodic and uniformly continuous in} 
	\quad \y \in \TdN.
\end{equation}
We can derive the above characterization  from the usual 
comparison  argument of viscosity solutions. (See, e.g. \cite{arisawaa}.) 
\end{remark}

\bigskip

\leftline{\bf Corollary 2.3}
\begin{corollary}
	Let $w_{\l}$ $(\l>0)$ be the solution of ($\ref{ergod}$)-($\ref{period}$) with
the Hamiltonian in Example 1.1 or 1.2. Then, for any $(x,p)\in \Omega\times \Rd$, ($\ref{wl}$) holds.
\end{corollary}

\bigskip
	
	In the rest part of this section, we shall prove Proposition 2.1, Theorem 2.2, 
Corollary 2.3 and Theorem A.
	We begin by recalling  the Jacobi's Lemma in the ergodic theory. (See e.g. \cite{arnold}.)

\bigskip

\leftline{{\bf Lemma 2.4}\quad (Jacobi)}

\begin{lemma}
	Let $\o=(\oa,\ob,...,\oN)\in \RN$ and consider the mapping $T$ from 
$\TN$ to $\TN$
\begin{equation}\label{map}
	T : \quad x \quad \to \quad \overline{x+\o} \equiv x+\o 		\quad 
	\pmod{1^N} \quad \in \TN.
\end{equation}
 Then, for each $x\in \TN$, $\{T^k x\}_{k\in \Nat}$ ($T^k$ denotes the $k$-th power of 
the mapping $T$) is dense in $\TN$ 
if and only if 
\begin{equation} \label{doku}
	\exists \quad \{z^n\}_{1\leq n\leq N}\in \ZN \quad
	\hbox{such that}\quad \Sigma_{n=1}^{N}z^n \on \in {\bf Z} \quad
	\hbox{implies}\quad  z^n=0\quad 1\leq \forall n\leq N.
\end{equation}
\end{lemma}

\bigskip

\leftline{\bf \it Proof of Proposition 2.1.}
We put for each $1\leq i\leq d$, $w_i=(\gia,...,\gin,...,\giN)\in \RN$, 
$$
	T_i \quad : \quad x\in \TN \quad \to \quad
	\overline{x+w_i} \equiv x+w_i \qquad \pmod{1^N}\quad \in \TN.
$$
Since the matrices $\Gn$ $(1\leq n\leq N)$ satisfy the Condition A, from Lemma 2.4
 we see that for any $x\in \TN$, $\{T_i^k x\}_{k\in {\bf N}}$ is dense in $\TN$ for each  $1\leq i\leq d$.\\

 Let us denote
$$
	\yo=(\yao,\ybo,...,\yNo),\qquad
	\yo'=({\overline{y}'}^1,{\overline{y}'}^2,...,{\overline{y}'}^N)\in \TdN,
$$
where
$$
	\yno=(\ynoa,\ynob,...,\ynod),\qquad 
	{\yo}'^n =({\yo}'^n_1,{\yo}'^n_2,...,{\yo}'^n_d)\in \Td \qquad 1\leq n\leq N,
$$
and  set for each $1\leq i\leq d$, 
$$
	p_i=(p_i^1,p_i^2,...,p_i^n,...,p_i^N) \qquad \qquad \qquad \qquad \qquad
	\qquad \qquad \qquad \qquad \qquad
$$
$$
	=(\g_i^1({\yo}'^1_i-\yo^1_i),\g_i^2({\yo}'^1_i-\yo^1_i)
	,...,\g_i^n({\yo}'^1_i-\yo^1_i),...,\g_i^N({\yo}'^1_i-\yo^1_i))
	\in \RN,
$$
$$\quad$$
$$
	q_i=(q_i^1,q_i^2,...,q_i^n,...,q_i^N) \qquad \qquad \qquad \qquad \qquad
	\qquad \qquad \qquad \qquad \qquad 
$$
$$
	=({\yo}'^1_i-\yo^1_i,{\yo}'^2_i-\yo^2_i
	,...,{\yo}'^n_i-\yo^n_i,...,{\yo}'^N_i-\yo^N_i)
	\in \RN,\qquad \qquad \quad
$$
$$\quad$$
$$
	\overline{p}_i\equiv p_i,\quad 
	\overline{q}_i\equiv q_i\qquad \pmod{1^N}\quad \in \TN.
$$
Then, for any $\d'>0$, from the above claim there exist 
$k_i=k(\d';p_i,q_i)$ $(1\leq i\leq d)$ such that
\begin{equation}\label{zzz}
	d_N(T^{k_i}\overline{p}_i,\overline{q}_i) \leq \d'.
\end{equation}
We set
\begin{equation}\label{defz}
	z_{\d'}=(\za,\zb,...,\zd)=(k_1,k_2,...,k_d)\in \Nat ^d\qquad 1\leq \forall i\leq d.
\end{equation}

\par

Next,  from the assumption of Proposition 2.1, there exists $T_{\d}>0$ such that
$$
	T_{\d}\leq \Delta_0(\d;\yo^1,{\yo'}^1+z_{\d'}),
$$
and a control $\a(\cdot)$ such that there exists a solution $\yaa(t)$ of
$$
	\frac{d}{dt}\yaa(t)=
	 b(x,\yaa(t),\Gb(\yaa(t)-\ya)+\yo^2,...,\GN(\yaa(t)-\ya)+\yo^N,\a(t))
	\quad t\geq 0,
$$
\begin{equation}\label{hora}
	\yaa(0)=\yo^1, \qquad |\yaa(T_{\d})-({\yo'}^{1}+z_{\d'})|\leq \d.
\end{equation}

Let $Y_{\a}(t)$ be the solution of ($\ref{system}$) in $\RdN$ with the same control $\a(\cdot)$ as in ($\ref{hora}$) with the initial condition
$$
	(\yaa(0),\yab(0),...,\yaN(0))=\yo.     
$$
From ($\ref{hora}$),  for each $1\leq n\leq N$
$$
	\yz^n(T_{\d}) - {\yo'}^n
	= \yo^n+ \Gn (\yz^1(T_{\d})-\yo^1)- {\yo}'^n
	\qquad \qquad \qquad \qquad \qquad
$$
$$
	= \yo^n+ \Gn ({\yo'}^1 + z_{\d} - \yo^1)- {\yo}'^n
	+\Gn (\yz^1(T_{\d})-({\yo'}^1 + z_{\d}))\qquad
$$
$$
	= (p^n_1,p^n_2,...,p^n_d)-(q^n_1,q^n_2,...,q^n_d) +\Gn z_{\d}
	+\Gn (\yz^1(T_{\d})-({\yo'}^1 + z_{\d}))
$$
$$
	=(p^n_1,p^n_2,...,p^n_d)+(k_1 \g_1^n,...,k_d \g_d^n)-(q^n_1,q^n_2,...,q^n_d)
	+\Gn (\yz^1(T_{\d})-({\yo'}^1 + z_{\d})).
$$
Therefore, from  $(\ref{zzz})$, $(\ref{hora})$, if we take $\d'=O(\d)$
$$
	d_{d}(\yz^n(T_{\d}), {\yo}'^n)\leq C\d,
$$
where $C>0$ is a constant independent on $\d$.
Therefore, $\yo$ is  approximately controllable to ${\yo}'$ in the system
$\overline{Y}_{\a}(t)$ with some estimate $\Delta(\d,\yo,{\yo}')$. The last inequality shows that $\Delta(\d,\yo,{\yo}')$  decreases to $0$ when $\Delta_0$ 
decreases to $0$.\\

\bigskip

\leftline{\bf \it Proof of Theorem 2.2}

(i) From Proposition 2.2, for any $\rho>0$, any $\yo, \quad \yo'$$\in \TdN$ , there exists  $0< \Delta(\d;\yo,\yo')< \rho$ such that  $\yo$ is approximately controllable to $\yo'$ with the estimate $\Delta(\d;\yo,\yo')$. Therefore, the condition (i) in Theorem A is satisfied 
and ($\ref{wl}$) holds.\\

(ii) Since the system ($\ref{first}$) is Lipschitz continuous (($\ref{Lip}$)), from ($\ref{zoo}$) $\overline{Y}_{\a}(t)$ is Lipschitz continuous. 
 Thus, from Proposition 2.2, any $\yo$ is approximately controllable to $\yo'$ with the estimate $\Delta(\d;\yo,\yo')$ in the Lipschitz continuous system $\overline{Y}_{\a}(t)$.
Therefore, the condition (ii) in Theorem A is satisfied 
and ($\ref{wl}$) holds.\\

(iii) From the continuity assumption for ($\ref{first}$), from ($\ref{zoo}$), 
$\overline{Y}_{\a}(t)$ is partially Lipschitz continuous. From the exact controllability assumption for the system ($\ref{first}$),  for any ${y'}^1\in \Td$, and for any $\yo$
$\in \TdN$, there exists a point ${\yo}'\in \TdN$ such that 
$$
	{\yo}'=({y'}^{1},...,{y'}^N),
$$
and $\yo$
is exactly controllable to ${\yo}'$.  Moreover, from Proposition 2.1, we have  the  approximated controllability of $\overline{Y}_{\a}(t)$.
Therefore,  the condition (iii)  in Theorem A is satisfied 
and ($\ref{wl}$) holds.\\

\leftline{\bf \it Proof of Corollary 2.3}

In Example 1.1, we rewrite ($\ref{qual}$) to
$$
	H(x,y,p)=\sup_{\a\in \Rd} \{-<\a,p\pd -\frac{1}{4a(x,y,\a)|\a|^2}-V(x,y)\},
$$
and see that $b(x,y,\a)=\a$ in ($\ref{first}$) defines a Lipschitz system  (uniformly in $(y^2,...,y^N)$). The function
$$
	h_{x,p}(y,\a)=<\a,p\pd + \frac{1}{4a(x,y,\a)|\a|^2}+V(x,y),
$$
satisfies ($\ref{bdd}$), ($\ref{power}$). 
For any $(y^2,...,y^N)$, any $\yo^1$ is approximately controllable to any ${\yo'}^1$ in the system ($\ref{first}$) with $b(x,y,\a)=\a$. Therefore, the condition (ii) in Theorem 2.2 
is satisfied, and we get ($\ref{wl}$).\\

In Example 1.2, we rewrite ($\ref{finite}$) to
$$
	H(x,y,p)=\sup_{|\a|\leq 1} \{-<a(x,y^1)\a,p\pd -V(x,y)\}.
$$
We see that the system ($\ref{first}$) with $b(x,y^1,y^2,...,y^N,\a)=a(x,y^1)\a$ satisfies the continuity assumption in (iii) in 
Theorem 2.2. The function
$$
	h_{x,p}(y,\a)=<a(x,y^1)\a,p\pd +V(x,y),
$$
satisfies ($\ref{bdd}$), ($\ref{power}$). Moreover, since $a(\cdot)>0$, for any $(y^2,...,y^N)$, there exists $T>0$ and any $\yo^1$ is exactly controllable to any ${\yo'}^1$ in the system ($\ref{first}$) with $b(x,y,\a)=a(x,y^1)\a$
within the time $T>0$. Therefore, the condition (iii) in Theorem 2.2 
is satisfied, and we get ($\ref{wl}$).\\

\bigskip

We shall give a brief  proof of Theorem A. We refer the reader to $\cite{arisawab}$ for more details.

\bigskip

\leftline{\bf \it Proof of Theorem A.}

We refer the reader to $\cite{arisawab}$ for the proof of (i). For (ii), the Lipschitz continuity of the system ($\ref{controllable}$)  leads 
$$
	|\l u_{\l}(y)-\l u_{\l}(y')|\leq C|y-y'|\qquad \forall y,\quad y'\in \TM,
$$
from the assumptions on $h$ (($\ref{bdd}$), ($\ref{power}$)). Therefore, we can extract a subsequence  $\l' u_{\l'}$  $(\l'\to 0)$ so that it converges uniformly to  a Lipschitz continuous function. On the other hand, from the approximated controllability, for any $y$, 
$y'\in \TM$, for any $\d>0$, there exist $T_{\d}>0$, and a control $\a(\cdot)$ such that
$d_M(y_{\a}(T_{\d}),y')\leq \d$
$$
	\l u_{\l}(y)\leq \l \int_{0}^{T_{\d}} h(y_{\a}(t),\a(t))dt 
	+ e^{-\l T_{\d}}\l u_{\l}(y_{\a}(T_{\d})).
$$
Since, $T_{\d}$ does not depend on $\l>0$, by combining the above two arguments we get 
$\lim_{\l \to 0}\l u_{\l}(y)=\hbox{Const.}$ with uniform convergence.  For (iii), the partial Lipschitz continuity of ($\ref{controllable}$)  leads 
$$
	|\l u_{\l}(y)-\l u_{\l}(y')|\leq C|y-y'|
$$
for any $y$, $y'\in \TM$ which have the  same  first $m$ entries. Since $y$ is exactly  controllable from $y$ to some $y'$, which has a specific first $m$ entries,  uniformly, 
$$
	\l u_{\l}(y)\leq \l \int_{0}^{T_0} h(y_{\a}(t),\a(t))dt 
	+ e^{-\l T_0}\l u_{\l}(y'),
$$
where $T_0<T$. By  combining the above arguments, we get 
$\lim_{\l \to 0}\l u_{\l}(y)=\hbox{Const.}$ with uniform convergence.\\
%%%%%%%%%%%%%%%%%%%%%%%%%%%%%%%%%%%%%%%%%%%%%%%%%%%%%%%%%%%%%ù

\newpage

\section{Cell problem in $\Rd$ and other related results}
	As we shall see in $\S$ 6, the formal asymptotic expansion for $\ue$ leads the   relationship ($\ref{ep}$).  We shall rigorously justify  the formal argument 
by the perturbed test function method in $\S$ 4. 
In the perturbed test function method, we need the so-called correctors as in 
\cite{evans}, \cite{lpv}.  We construct the approximated correctors in ($\ref{pp}$)-($\ref{ll}$) from the functions $\wu$, $\wo$ in Remark 2.2.\\

\bigskip

\leftline{\bf Theorem 3.1}

\begin{proposition}
	 Let $\wu$, $\wo$ be the solutions of ($\ref{ari}$), ($\ref{sawa}$) respectively, 
satisfying ($\ref{ppp}$).
Then, the functions $\vu$, $\vo$ defined by
\begin{equation}\label{vu}
	\vu (\ya)= \wu(\Ga \ya, \Gb \ya,...,\GN \ya)\qquad \ya \in \Rd,
\end{equation}
\begin{equation}\label{vo}
	\vo (\ya)= \wo(\Ga \ya, \Gb \ya,...,\GN \ya)\qquad \ya \in \Rd,
\end{equation}
are respectively  a viscosity subsolution  and a viscosity supersolution of 
($\ref{pp}$), ($\ref{ll}$).
\end{proposition}

\bigskip

\leftline{\bf \it Proof of Theorem 3.1}

	We shall prove that $\vu (\ya)$ is a viscosity subsolution of ($\ref{pp}$).  Let $\phi (\ya)\in C^2(\Rd)$ and assume  that  $(\vu -\phi) (\ya)$
takes a local strict maximum at ${\yo}^{1}\in \Rd$ in a neighborhood 
$U$ of ${\yo}^{1}$. We are to show that
\begin{equation}\label{neig}
	-\overline{H}(x,p)+ H(x,\Ga {\yo}^{1},\Gb {\yo}^{1},...,\GN {\yo}^{1},p+\Na \phi 	(\overline{y}^1)) \leq \d.
\end{equation}
	For any $\beta >0$, set
$$
	{\Psi}_{\beta}(\ya,\yb,...,\yN)= 
	\wu (\ya,\yb,...,\yN)-\phi(\ya)
	-\beta  \sum_{i=1,n=2}^{i=d,n=N} |y_i^n - \gin y_i^1|^2
$$
$$
	\qquad \qquad \qquad \qquad \qquad \qquad \qquad \qquad \qquad 
	\quad \hbox{for} \quad (\ya,\yb,...,\yN)\in \RdN.
$$
Let ${\Psi}_{\beta}$ take its strict maximum in the compact set
$$
	U'=\{y=(\ya,\yb,...,\yN)\quad|\quad
	|(\ya,\yb,...,\yN)-(\Ga \overline{y}^1,\Gb \overline{y}^1
	,...,\GN \overline{y}^1)|\leq \delta  	\}\subset \RdN
$$ 
at $(y^1_{\beta},y^2_{\beta},...,y^N_{\beta})\in U'$.
Then, by denoting $y^n_{\beta},=(y^n_{\beta-1},y^n_{\beta_2},...,y^n_{\beta_d})$ $(1\leq n\leq N)$ we have

$$
	\wu (y^1_{\beta},y^2_{\beta},...,y^N_{\beta})-\phi(y^1_{\beta})
	-\beta \Sigma_{i=1,n=2}^{i=d,n=N} |y^n_{\beta_i}-\gin y^1_{\beta_i}|^2
	\qquad \qquad \qquad \qquad \qquad \qquad
$$
$$\quad$$
$$
	\geq 
	\wu (\Ga\overline{y}^1,\Gb\overline{y}^1,...,\GN\overline{y}^1)
	-\phi(\Ga\overline{y}^1)
	-\beta \Sigma_{i=1,n=2}^{i=d,n=N} |\gin \overline{y}^1_i
			-\gin \g^1_i\overline{y}^1_i|^2
	\qquad \qquad
$$
$$\quad$$
$$
	\geq
	\wu (\Ga\overline{y}^1,\Gb\overline{y}^1,...,\GN\overline{y}^1)
	-\phi((\Ga)^{-1}\Ga\overline{y}^1),
	\qquad \qquad \qquad \qquad \qquad \qquad \qquad
$$
$$\quad$$
and we have
$$
	\beta \Sigma_{i=1,n=2}^{i=d,n=N} |y^n_{\beta_i}-\gin y^1_{\beta_i}|^2
	\leq M \qquad \hbox{for}\quad \forall \beta>0,
$$
where $M>0$ is a constant. 
Thus, we get 
\begin{equation}\label{beta}
	\lim_{\beta \to \infty}(y^1_{\beta},y^2_{\beta},...,y^N_{\beta})\to 
	(\Ga\overline{y}^1,\Gb\overline{y}^1,...,\GN\overline{y}^1).
\end{equation}

Since $\wu$ is the viscosity subsolution of $(\ref{ari})$, by setting
$$
	\phi_{\beta}(\ya,\yb,...,\yN)
	=\phi(\ya)+\beta
	 \Sigma_{i=1,n=2}^{i=d,n=N} |\yy - \gin y^1_i|^2,
$$
from the definition of viscosity solutions 
\begin{equation}\label{wuvis}
	- \overline{H}(x,p) + 
	H(x,y^1_{\beta},y^2_{\beta},...,y^N_{\beta},p+ \Sigma_{n=1}^{N} \Gn\Nn\phi_{\beta} 		(y^1_{\beta},y^2_{\beta},...,y^N_{\beta}))\leq \d.
\end{equation}

For $n=1$, we have
$$
	\Ga \N_{y_i^1} \phi_{\beta}(\ya,\yb,...,\yN)
	=\N_{y_i^1} \phi (\ya)
	-2\beta \Sigma_{m=2}^{N} \g_i^m(y_i^m-\g_i^m y^1_{i}),
$$
and for $n\neq 1$, we have
$$
	\Gn \Nin \phi_{\beta}(\ya,\yb,...,\yN)
	=2\beta \gin(\yy-\gin y_i^1).
$$
By introducing the above into $(\ref{wuvis})$ and by letting $\beta \to \infty$,
from ($\ref{beta}$) we get
$$
	- \overline{H}(x,p)+
	H(x,\Ga {\yo}^{1},\Gb {\yo}^{1},...,\GN {\yo}^{1},p+\Na \phi (\overline{y}^1),\a) 
	\leq \d.
$$

The proof that $\vo$ is the viscosity of $(\ref{ll})$ is similar to above, and 
we do not write it here.\\

\bigskip

\bigskip

	We collect  the following  results similar to those in  \cite{evans} and \cite{lpv}. We shall use some of them  in $\S\S$ 4, 5.
	
\bigskip

\leftline{\bf Proposition 3.2}
\begin{proposition}
Consider ($\ref{ergod}$)-($\ref{period}$), and assume that ($\ref{wl}$) holds. Then, $\Ho(x,p)$ defined by ($\ref{ep}$) has  following properties.\\

(i) If $H(\cdot)$ is locally Lipschitz:
$$
	|H(x,y,p)-H(x',y',p')|\leq C(L') |(x,y,p)-(x',y',p')|
$$
$$
	\forall (x,y,p),\quad (x',y',p')\in \{|x|\leq L'\}\times \RdN \times \{|p|\leq L'\},
$$
for $\forall L'>0$, then so is $\Ho(\cdot)$:
$$
	|\Ho(x,p)-\Ho(x',p')|\leq \overline{C}(L')|(x,p)-(x',p')|
$$
$$
	\forall (x,p),\quad (x',p')\in \{|x|\leq L'\}\times \{|p|\leq L'\},
$$ 
for  $\forall L'>0$, where the Lipschitz constant  $\overline{C}(L')$ is determined only by $C(L')$.\\

(ii) Assume that for any  $\d>0$, there exists a  Lipschitz continuous subsolution $\wu$ of ($\ref{ari}$). Then, if for any $(x,\ya,...,\yN)\in \Omega \times \RdN$ $H$ is convex in $p\in \Rd$, then for any $x\in \Rd$ $\Ho$ is convex in $p\in \Rd$.\\

(iii) If there exists $\theta\geq 0$ such that 
$$
	\lim_{|p|\to \infty}|p|^{-\theta}H(x,y^1,...,y^N,p)= \infty \quad
$$
$$
	 \hbox{uniformly in} \quad (x,y^1,...,y^N)\in \Omega \times \{|y|\leq L',\quad y\in 	 \RdN\} \quad \hbox{for any}\quad L'>0,
$$
then
$$
	\lim_{|p|\to \infty}|p|^{-\theta}\Ho (p,x)= \infty \quad \hbox{uniformly in} \quad
	x\in \Omega. 
$$

\end{proposition}

\bigskip

\leftline{\bf \it Proof of Proposition 3.2.}

1. The proof is  similar to that of \cite{evans} and \cite{lpv}, so we do not repeat it here.\\

2. Let $x\in \Omega$, $p,p'\in \Rd$. From the assumption and Remark 2.2, for any $\d>0$  there exist 
$w^p$, $w^{p'}$, $w^{\frac{p+p'}{2}}$ Lipschitz in $\y\in \RdN$ such that
$$
	-\Ho(x,p)+  H(x,\ya,\yb,...,\yN,p+ \Sigma_{n=1}^N \Gn \Nn w^p (y))
	 \leq \d
	\qquad y=\y \in \RdN,
$$
$$\quad$$
$$
	-\Ho(x,p')+  H(x,\ya,\yb,...,\yN,p'+ \Sigma_{n=1}^N \Gn \Nn w^{p'} (y))
	 \leq \d
	\qquad y=\y \in \RdN,
$$
$$\quad$$
where $w^p$, $w^{p'}$  are periodic in $\y\in \TdN$.
We are to show that
\begin{equation}\label{convex}
	\Ho(x,\frac{p+p'}{2})\leq \frac{\Ho(x,p)+\Ho(x,p')}{2}.
\end{equation}
For this purpose, we assume that
\begin{equation}\label{apple}
	\Ho(x,\frac{p+p'}{2})>     \frac{\Ho(x,p)+\Ho(x,p')}{2},
\end{equation}
and we shall look for a contradiction.
Let $\rho$ be a mollifier ($\rho \in C^{\infty}_0(\RdN)$, $\hbox{supp}\rho \subset \subset
[-1,1]^{d\times N}$, $\int_{-\infty}^{\infty} \rho (y')dy'=1$,)
and set
$$
	w_{\delta}(y)=\int_{-\infty}^{\infty}
	\rho(\frac{y-y'}{\d}) \frac{w^p+w^{p'}}{2}(y')dy' \qquad y\in \RdN.
$$
 Then,
$$
	H(x,\ya,\yb,...,\yN, \frac{p+p'}{2}+\Sigma_{n=1}^N \Gn \Nn w_{\delta}(y))
	\leq 
$$
$$
	\leq \int_{U_{\d}(y)}\rho(\frac{y-y'}{\d})
	H(x,{\ya},{\yb},...,{\yN},
	 \frac{p+p'}{2}+\Sigma_{n=1}^N \Gn \Nn \frac{w^p+w^{p'}}{2}(y'))dy'
$$
$$
	= \int_{U_{\d}(y)}\rho(\frac{y-y'}{\d})
	H(x,{\ya}',{\yb}',...,{\yN}',
	 \frac{p+p'}{2}+\Sigma_{n=1}^N \Gn \Nn \frac{w^p+w^{p'}}{2}(y'))dy'+o(1)
$$
$$
	\leq \frac{1}{2}\int_{U_{\d}(y)}\rho(\frac{y-y'}{\d})
	H(x,{\ya}',{\yb}',...,{\yN}',
	 p+\Sigma_{n=1}^N \Gn \Nn w^p(y'))dy'\qquad \qquad
$$
$$
		+\frac{1}{2}\int_{U_{\d}(y)}\rho(\frac{y-y'}{\d})
	H(x,{\ya}',{\yb}',...,{\yN}',
	 p'+\Sigma_{n=1}^N \Gn \Nn w^{p'}(y'))dy'+o(1)
$$
$$
	=\frac{1}{2}\Ho(x,p)+\frac{1}{2}\Ho(x,p')+o(1)+2\d,
$$
since $w^p$, $w^{p'}$ are Lipschitz. Now, by letting $\d \to 0$  
from $(\ref{apple})$
$$
	H(x,\ya,\yb,...,\yN,\frac{p+p'}{2}+\Sigma_{n=1}^N \Gn \Nn \frac{w^p+w^{p'}}{2}(y))
$$
$$
	\leq \frac{1}{2}\Ho(x,p)+\frac{1}{2}\Ho(x,p')
	< \Ho(x,\frac{p+p'}{2}) \quad \forall y=\y\in \RdN.
$$
However, since $\Ho(x,\frac{p+p'}{2})$ is the least value of $r\in \R$ such that the 
following problem with $(\ref{period})$ admits a viscosity subsolution (see \cite{infinite}, \cite{apl})
$$
	H(x,\ya,\yb,...,\yN,\frac{p+p'}{2}+\Sigma_{n=1}^N \Gn \Nn w(y))
	\leq r,
$$
 we get a contradiction. Thus, we have proved our claim.
\newline
3.  The proof is same to \cite{evans} and \cite{lpv},  and we do not write it here.\\

\section{Hologenization problem.}

The following is a main result of this paper. 
\leftline{Theorem 4.1}
\begin{theorem}
Assume that $b(\cdot)$, $g(\cdot)$ in (1), (2) satisfy either one of the conditions in Theorem 2.2. Then as $\e$ goes to $0$, the solutions $u_{\e}(x,t)$, $u_{\e}(x)$ 
of (1), (2) converge respectively to the solutions $\overline{u}(x,t)$, $\overline{u}(x)$ of (7), (8) uniformly in $\Omega\times [0,T]$, $\Omega$ for any $T>0$. The effective Hamiltonian $\overline{H}$ in (7), (8) is defined by (11).  
\end{theorem}

\leftline{\bf Remark 4.1}
\begin{theorem}
The effective Hamiltonian $\overline{H}(\cdot)$ satisfies the properties in Proposition 3.2. 
\end{theorem}

\leftline{\bf Corollary 4.2}
\begin{theorem}
The statements in Examples 1.1 and 1.2 hold. In Example 1.1, $\overline{H}$ is locally Lipschitz continuous, convex in $p$ and satisfies 
$$
	\lim_{|p|\to \infty} |p|^{-2} \overline{H}(x,p)=\infty \quad\hbox{uniformly in}\quad x\in \Omega.
$$
In Example 1.2, $\overline{H}$ is locally Lipschitz continuous, convex in $p$ and satisfies 
$$
	\lim_{|p|\to \infty} \overline{H}(x,p)=\infty \quad\hbox{uniformly in}\quad x\in \Omega.
$$
\end{theorem}
Now, we enter into the proofs. \\

$Proof$ $of$ $Theorem$ $4.1$ 
  From the assumptions for ${H}$, there exist subsequences $\{u_{\e'}(x,t)\}$, $\{u_{\e'}(x)\}$ which converge to $\overline{u}(x,t)$, $\overline{u}(x)$ uniformly in $\Omega\times [0,T]$, $\Omega$, as $\e'\to 0$ for any $T>0$. We are to show that $\overline{u}(x,t)$, $\overline{u}(x)$ satisfy (7), (8) with Hamiltonian $\overline{H}(\cdot)$ which is well-definded by Theorem 2.2. This also proves the uniqueness of limits $\overline{u}(x,t)$, $\overline{u}(x)$. We shall show the proof only for the stationary case (the evolutionary case is proved similarly), that $\overline{u}(x)$ is a subsolution of (8). The proof that $\overline{u}(x)$ is a supersolution is similar, and we do not write it here. \\

Let $\phi(x)\in C^2(\Omega)$ and assume that $(\overline{u}-\phi)(x)$ takes a local strict maximum at $x_0\in \Omega$. We shall show 
\begin{equation}
	\mu \phi(x_0) + \overline{H}( x_0,\n \phi(x_0))\leq 0, \quad x\in \Omega. 
\end{equation}
For this purpose, let us assume that (68) does not hold: 
\begin{equation}
	\mu \phi(x_0) + \overline{H}( x_0,\n \phi(x_0))=3 \d > 0,
\end{equation}
and we shall look for a contradiction. There exists a neighborhood $U$ of $x_0$ such that 
$$
	\mu \phi(x) + \overline{H}( x,\n \phi(x))\geq  \d \quad \forall x\in U.
$$
Let $p_0=\n \phi(x_0)$. For $(x_0,p_0)$, from Theorem 3.1 we can take $\overline{v}(y^1)$ a viscosity supersolution of (15): 
\begin{equation}
	H(x_0,\Gamma^1 \ya,\Gamma^2 \ya,...,\Gamma^N \ya,p_0+\n_{y^1}\overline{v}(y^1))\geq \overline{H}( x_0,p_0)-\delta \quad y^1\in {\bf R^d}. 
\end{equation}
Set 
$$
	\phi_{\e}(x)=\phi(x)+ \e^{1}_{1} \overline{v}(\frac{x}{\e^1})= \phi(x)+ \e^{1}_{1} \overline{v}(\frac{x_1}{\e^1_1}, \frac{x_2}{\e^1_1},...,\frac{x_d}{\e^1_1} )
	\quad (x\in \overline{\Omega}). 
$$
We claim that $\phi_{\e}(x)$ defined above satisfies 
$$
	\mu \phi_{\e}(x)+ H(x,\frac{x}{\e^1},...,\frac{x}{\e^n},...\frac{x}{\e^N}, \n \phi_{\e}(x)) \geq \delta \quad x\in U_r(x_0), 
$$
in the sense of viscosity solutions if we take $r>0$ small enough. Therefore, we have to show that there exists $r>0$ such that for any $\psi(x)$$\in C^2({\bf R^d})$
 such that $\phi_{\e}-\psi$ takes its local minimum at $\overline{x}\in U_r(x_0)$ (we may assume that $\phi_{\e}(\overline{x})=\psi(\overline{x})$), the following holds. 
$$
	\mu \psi(\overline{x})+ H(\overline{x},\frac{\overline{x}}{\e^1},...,\frac{\overline{x}}{\e^n},...\frac{\overline{x}}{\e^N}, \n \psi(\overline{x})) \geq \delta. 
$$
If 
$$
	(\phi_{\e}-\psi)(x)=\e^{1}_{1} \overline{v}(\frac{x}{\e^1_1}) - (\psi(x)-\phi(x)) 
$$
takes its local minimum at $x=\overline{x}$ in $U_s(x_0)$, then by putting $y^1=\frac{x}{\e^1}$
$$
	\eta(y^1)=\frac{1}{\e^{1}_{1}}(\psi-\phi)(\e^{1}_{1}y^1_1,\e^{1}_{1}y^1_2,...,\e^{1}_{1}y^1_d ), 
$$
$(\overline{v}-\eta)(y^1)$ takes a local minimum at 

$$
	\overline{y}^1=\frac{\overline{x}}{\e^{1}}= (\frac{\overline{x_1}}{\e^{1}_1},\frac{\overline{x_2}}{\e^{1}_1},...\frac{\overline{x_d}}{\e^{1}_1} ). 
$$
Since $\overline{v}(y^1)$ is the viscosity supersolution of (70), and since 
$$
	\n_{y^1}\eta(\overline{y}^1)=\n (\psi-\phi)(\e^1_1\Gamma^1\overline{y}^1) = \n (\psi-\phi)(\overline{x})
$$
we have 
$$
	-\overline{H}(x_0,p_0)+ H(x_0,\Gamma^1 \overline{y}^1,\Gamma^2 \overline{y}^1,...,\Gamma^N \overline{y}^1,p_0+\n_{y^1} \eta(\overline{y}^1))
$$
$$
	= -\overline{H}(x_0,p_0)+ H(x_0,\frac{\overline{x}}{\e^{1}},\frac{\overline{x}}{\e^{2}},...,\frac{\overline{x}}{\e^{N}},p_0+\n(\psi-\phi)(\overline{x}) )\geq -\d. 
$$
Hence, from (69), 
$$
	\mu\phi(x_0)+ H(x_0,\frac{\overline{x}}{\e^{1}},\frac{\overline{x}}{\e^{2}},...,\frac{\overline{x}}{\e^{N}},p_0+\n(\psi-\phi)(\overline{x}) )\geq 2\delta. 
$$
Hence, if we take $s>0$ small enough we get (71). Now, from the comparison result of the viscosity solutions of (2) and (71), we have 
$$
	(u_{\e}-\phi_{\e})(x_0) \leq \max_{\p U_r(x_0)} (u_{\e}-\phi_{\e}) + \delta, 
$$
and by letting $\e\to 0$, we have 
$$
	(\overline{u}-\phi)(x_0) \leq \max_{\p U_r(x_0)} (\overline{u}-\phi) + \delta. 
$$
Since $\delta>0$ is arbitrary, 
$$
	(\overline{u}-\phi)(x_0) \leq \max_{\p U_r(x_0)} (\overline{u}-\phi), 
$$
which is a contradiction to the assumption that $x_0$ is a local strict maximum of $\overline{u}-\phi$. Thus, we have proved that (69) leads a contradiction and that 
$\overline{u}$ is a viscosity subsolution of (8).\\

$Proof$ $of$ $Corollary$ $4.2$\\
From Corollary 2.3, the effective Hamiltonians are well-definded. Thus, by using Theorem 4.1, we confirm the statements in Examples 1.1 and 1.2 are true. Moreover, from Proposition 3.2, the effective Hamiltonians have the properties in the Corollary. \\

\section{Quasi periodic and almost periodic homogenizations.}
First, we study (20), (21) with the quasi periodic Hamiltonian $F(x,y,p)$ in (22). We set 
$$
	\gamma^n_i=\frac{T^n_i}{T^1_i} \quad 1\leq i\leq d,\quad 1\leq n\leq N, 
$$
and define $d\times d$ matrices $\Gamma^n$ ($1\leq n\leq N$) by using the above $\gamma^n_i$ in (10). We also set 
$$
	b'(x,y^1,y^2,...,y^N,\a)= \sum_{n=1}^{N} b^n(x,(\Gamma^n)^{-1}y^n,\a), 
$$
$$
	g'(x,y^1,y^2,...,y^N,\a)= \sum_{n=1}^{N} g^n(x,(\Gamma^n)^{-1}y^n,\a), 
$$
and define the Hamiltonian $H'$ on $\Omega\times {\bf R^{d\times N}}\times {\bf R^d}$ by (3) with $b=b'$, $g=g'$. (Remark that $H'$ satisfies (4), (5), and (6).)\\
Our result for the quasi periodic Homogenization is the following. 

\leftline{\bf Theorem 5.1}
\begin{theorem}
Consider (20), (21) with F in (22). Assume that $T^1_i=1$ ($1\leq \forall i\leq d$), and that matrices $\Gamma^n$ ($1\leq n\leq N$) satisfy the Condition A. Then, the solutions $u_{\e}(x,t)$,  $u_{\e}(x)$ converge respectively to $\overline{u}(x,t)$, $\overline{u}(x)$ of (30), (31) with the effective Hamiltonian $\overline{F}$ given by (32) uniformly on $\Omega\times [0,T]$, $\Omega$ for any $T>0$. The effective Hamiltonian $\overline{F}$ is characterized by (34). 
\end{theorem}

Next, we consider (20), (21) with the class $\mathcal{B}_0$ Hamiltonian $F$ in (25). We denote by $\overline{F}^N$ the effective Hamiltonian for ${F}^N(x,y,p)$ (see Theorem 5.1). Our homogenization result for the class $\mathcal{B}_0$ Hamiltonian is the following.

\leftline{\bf Theorem 5.2}
\begin{theorem}
Consider (20), (21) with F in (25). The solutions of (20), (21), $u_{\e}(x,t)$,  $u_{\e}(x)$ converge respectively to $\overline{u}(x,t)$, $\overline{u}(x)$ of (30), (31) with the effective Hamiltonian $\overline{F}$ given by (35) uniformly on $\Omega\times [0,T]$, $\Omega$ for any $T>0$. 
\end{theorem}

Finally, we study (20), (21) with the class $\mathcal{B}_1$ Hamiltonian $F$ in (26). Our homogenization result for the class $\mathcal{B}_1$ Hamiltonian is the following. 

\leftline{\bf Theorem 5.3}
\begin{theorem}
Consider (20), (21) with F in (26). The solutions of (20), (21), $u_{\e}(x,t)$,  $u_{\e}(x)$ converge respectively to $\overline{u}(x,t)$, $\overline{u}(x)$ of (30), (31) with the effective Hamiltonian $\overline{F}$ given by (36) uniformly on $\Omega\times [0,T]$, $\Omega$ for any $T>0$. 
\end{theorem}

\leftline{\bf Corollary 5.4}
\begin{theorem}
The statements in Examples 1.3 and 1.4 hold. 
\end{theorem}

Now, we enter into the proofs of above theorems.\\

$Proof$ $of$ $Theorem$ $5.1.$ \\
First, we remark that there exists a unique solution $v_{\l}$ of (33) which grows at most linearly and in $UC({\bf R^d})$. From the assumption on $F$: (6), 
$|\n_y v_{\l}|_{\infty}<C$ uniformly in $y\in {\bf R^d}$. Therefore, there exists a subsequence $\l'\to 0$ such that 
$$
	\lim_{\l'\to 0} \l' v_{\l'}(y)=\exists d_{x,p}\quad \hbox{locally uniformly in}\quad y\in {\bf R^d}. 
$$
By using a similar argument as in \S 3, we see the relationship between the solution $w_{\l}(y^1,...,y^N)$ of 
$$
	\l w_{\l} + H'(x,y^1,y^2,...,y^N,p+\sum_{n=1}^{N}  \Gamma^n \n_{y^n} w_{\l})=0\in {\bf T^{d\times N}}
	\quad (y^1,...,y^N)\in {\bf R^{d\times N}}, 
$$
and the solution $v_{\l}(y)$ of (33), that is 
$$
	v_{\l}(y)=w_{\l}(y,\Gamma^2 y,...,\Gamma^N y) \quad y\in {\bf R^d}. 
$$
Therefore, by using the fact that $w_{\l}$ is periodic, we see that $v_{\l}$ is bounded. Moreover, from Theorems 2.2, 3.1, 4.1, we see the uniform convergence of 
$\l' v_{\l'}$ to $d_{x,p}$ and (34) by taking $\underline{v}(y)$, $\overline{v}(y)$$=v_{\l'}(y)$ for $\l'>0$ small enough. We can prove the uniqueness of $d_{x,p}$ by using the usual comparison argument and (34). (See \cite{arisawaa} for such a characterization.) Therefore, we can use $\underline{v}(y)$, $\overline{v}(y)$ as approximated correctors in the perturbed test function method as in \S 4, and can show that the limits $\overline{u}(x,t)$, $\overline{u}(x)$ satisfy (30), (31) respectively.\\

$Proof$ $of$ $Theorem$ $5.2.$\\
First, from the assumption on $F$, $\overline{F}^N(x,p)$ $(N\in {\bf N})$ is uniformly bounded and we can extract a subsequence $\overline{F}^{N'}(x,p)$ which converges to $\overline{F}(x,p)$. To see the uniqueness of $\overline{F}(x,p)$, assume that there exist two subsequences $\overline{F}^{N'}(x,p)$, $\overline{F}^{N''}(x,p)$ which converge to $\overline{F}'(x,p)$, $\overline{F}''(x,p)$ respectively. That is for any $\delta>0$, if we take $N'$, $N''$ large enough, from Theorem 5.1, there exist uniformly continuous bounded functions  $\underline{v}^{N'}(y)$, $\overline{v}^{N'}(y)$,  $\underline{v}^{N''}(y)$, $\overline{v}^{N''}(y)$ which satisfy 
$$
	F^{N'}(x,y,p+\n_y \underline{v}^{N'}(y))\leq \overline{F}'(x,p) +\delta, 
$$
$$
	F^{N'}(x,y,p+\n_y \overline{v}^{N'}(y))\geq \overline{F}'(x,p) -\delta, 
$$
$$
	F^{N''}(x,y,p+\n_y \underline{v}^{N''}(y))\leq \overline{F}''(x,p) +\delta, 
$$
$$
	F^{N''}(x,y,p+\n_y \overline{v}^{N''}(y))\leq \overline{F}''(x,p) -\delta. 
$$
Since (6) holds uniformly in $N\in {\bf N}$, we see that 
$$
	|\n_y \underline{v}^{N'}|, \quad |\n_y \overline{v}^{N'}|, \quad |\n_y \underline{v}^{N''}|, \quad |\n_y \overline{v}^{N''}| \leq C. 
$$
Therefore, from the definition of $F$, we have 
$$
	F(x,y,p+\n_y \underline{v}^{N'}(y))\leq \overline{F}'(x,p) +2\delta, 
$$
$$
	F(x,y,p+\n_y \overline{v}^{N'}(y))\geq \overline{F}'(x,p) -2\delta, 
$$
$$
	F(x,y,p+\n_y \underline{v}^{N''}(y))\leq \overline{F}''(x,p) +2\delta, 
$$
$$
	F(x,y,p+\n_y \overline{v}^{N''}(y))\leq \overline{F}''(x,p) -2\delta. 
$$
Since  $\underline{v}^{N'}(y)$, $\overline{v}^{N'}(y)$,  $\underline{v}^{N''}(y)$, $\overline{v}^{N''}(y)$ are bounded in $y\in {\bf R^d}$, by using the comparison result in the unbounded domain (see \cite{crandall}), we have $\overline{F}'(x,p)$$=\overline{F}''(x,p)$. \\

Next, from the assumption (6), as usual, there exist subsequences of $u_{\e}(x,t)$, $u_{\e}(x)$ such that $\lim_{\e'\to 0}u_{\e'}(x,t)$$=\overline{u}(x,t)$,  
 $\lim_{\e'\to 0}u_{\e'}(x)$$=\overline{u}(x)$ uniformly on $\Omega\times [0,T]$, $\Omega$ for any $T>0$. We are to show that $\overline{u}(x,t)$, $\overline{u}(x)$
 are the viscosity solution of (30), (31), which also implies the uniqueness of the limits $\overline{u}(x,t)$, $\overline{u}(x)$. In the following, we shall only prove the stationary case, the evolutionary case can be proved in a similar way. Let $\phi(x)$$\in C^2({\bf R^d})$ and assume that $(\overline{u}-\phi)(x)$ takes a strict local maximum at $x_0\in \Omega$. We are to show that  
\begin{equation}
	\mu\phi(x_0)+  \overline{F}(x_0,\n \phi(x_0))\leq 0. 
\end{equation}
For this purpose, we assume that (72) does not hold and we shall look for a contradiction. Hence, assume that 
\begin{equation}
	\mu\phi(x_0)+  \overline{F}(x_0,\n \phi(x_0))=3\delta>  0. 
\end{equation}
There exists a neighborhood $U$ of $x_0$ such that 
$$
	\mu\phi(x)+  \overline{F}(x,\n \phi(x))\geq 2\delta, \quad \forall x\in U.
$$
For the above $x_0$, put $p_0=\n \phi(x_0)$. From the definition of $\overline{F}(x_0,p_0)$ and from Theorem 5.1, if we take $N$ large enough there exists $\overline{v}^N(y)$ a bounded uniformly continuous viscosity solution of 
$$
	-\overline{F}(x_0,p_0)+ F^N(x_0,y,p_0+\n_{y} \overline{v}^N(y)) \geq -\delta \quad y\in {\bf R^d}. 
$$
Again, by using the definition of $F$, $\overline{v}^N(y)$  satisfies 
\begin{equation}
	-\overline{F}(x_0,p_0)+ F(x_0,y,p_0+\n_{y} \overline{v}^N(y)) \geq -2 \delta \quad y\in {\bf R^d}.
\end{equation}
Now, set 
$$
	\phi_{\e}(x)=\phi(x)+ \e^1_1 \overline{v}^N (\frac{x}{\e^1}) \quad x\in \overline{\Omega},
$$
and we claim that $\phi_{\e}(x)$ satisfies 
\begin{equation}
	\mu\phi_{\e}(x)+ F(x,\frac{x}{\e},\n \phi_{\e}(x)) \geq \delta \quad x\in U_s(x_0),
\end{equation}
if we take $s>0$ small enough, in the sense of viscosity solutions. The inequality (75) is confirmed similarly as in the argument in the proof of Theorem 4.1, and we do not repeat it here. Thus, from (21) and (75) by using the comparison result, 
$$
	(u_{\e}-\phi_{\e})(x_0) \leq \max_{\p U_s(x_0)} (u_{\e}-\phi_{\e}) + \delta. 
$$
and by letting $\e\to 0$, we have 
$$
	(\overline{u}-\phi)(x_0) \leq \max_{\p U_s(x_0)} (\overline{u}-\phi) + \delta, 
$$
where $\d>0$ is arbitrary. Thus, we get a contradiction to the assumption that $x_0$ is a strict maximum of $\overline{u}-\phi$. Therefore $\overline{u}$ is a viscosity subsolution of (31). The proof that $\overline{u}$ is a viscosity supersolution of (31) is similar, and we do not write it here. \\

$Proof$ $of$ $Theorem$ $5.3.$\\
As in the proof of Theorem 5.1, there exists a subsequence $\l'\to 0$ such that 
$$
	\lim_{\l'\to 0} \l' v_{\l'}(y)=\exists d_{x,p}\quad \hbox{locally uniformly in}\quad y\in {\bf R^d}, 
$$

$$
	\lim_{\l'\to 0} (v_{\l'}(y)-v_{\l'}(0))=\exists v(y) \quad \hbox{locally uniformly in}\quad y\in {\bf R^d}. 
$$
We have 
$$
	d_{x,p}+ F(x,y,p+\n_y v(y))=0\quad y\in {\bf R^d},
$$
where $v$ is uniformly continuous and grows at most linearly. We shall see that $v$ satisfies (37). \\

Let $\overline{v}_{\l}(y)$ be the solution of 
$$
	\l \overline{v}_{\l}(y)+\sup_{\a\in \mathcal{A}} \{
	\left\langle -\a,\n_y \overline{v}_{\l}(y) \right\rangle_{{\bf R^d}\times {\bf R^d}} -V(x,y,\a) - \left\langle -\a,p\right\rangle_{{\bf R^d}\times {\bf R^d}}
	+c(x,p)
	\}=0,
$$
which grows at most linearly at infinity. Obviosly, 
\begin{equation}
	v_{\l}(y)=\overline{v}_{\l}(y)-\frac{c(x,p)}{\l}, 
\end{equation}
and for any $y$, $y'$$\in {\bf R^d}$, by putting $z=\frac{y-y'}{y-y'}$, we have 
$$
	\overline{v}_{\l}(y) \leq \int_{0}^{|y-y'|} e^{-\l t} \{
	V(x,y+tz,z) + \left\langle z,p\right\rangle_{{\bf R^d}\times {\bf R^d}}
	-c(x,p)
	\}dt + e^{-\l |y-y'|}\overline{v}_{\l}(y').
$$
Thus, by letting $\l'\to 0$, for any $y\in {\bf R^d}$ $(y\neq 0)$
$$
	\frac{|v(y)|}{|y|}=\lim_{\l' \to 0} \frac{v_{\l'}(y)-v_{\l'}(0)}{|y|} = \lim_{\l' \to 0} \frac{\overline{v}_{\l'}(y)-\overline{v}_{\l'}(0)}{|y|}
$$
$$
	\leq \frac{1}{|y|} \int_{0}^{|y|}  \{
	V(x,y+t\frac{y}{|y|},\frac{y}{|y|}) + \left\langle \frac{y}{|y|},p\right\rangle_{{\bf R^d}\times {\bf R^d}}
	-c(x,p)
	\}dt \leq 0, 
$$
where the last inequality holds from (27), (28). Hence, we have proved (37). \\

Since $v$ satisfies (37), we can use $v$ as a corrector in the perturbed test function method as in \S 4, and we can show that the limit $\overline{u}(x,t)$, $\overline{u}(x)$ satisfy (30), (31) respectively with the effective Hamiltonian in (36). \\

$Proof$ $of$ $Corollary$ $5.4.$\\
Since the Hamiltonians $F$ in Examples 1.3, 1.4 are in class $\mathcal{B}_1$, we can apply the result in Theorem 5.3. \\

 \leftline{\bf Remark 5.2}
\begin{remark}
The uniqueness of the effective Hamiltonian $\overline{F}$ in Theorem 5.3 is an open problem. (Remark that in Theorems 4.1, 5.1, and 5.2, we proved the uniqueness of effective Hamiltonians.) The difficulty comes from the unboundedness of the corrector $v$ in (34). 
\end{remark}

\section{Formal asymptotic expansions.}

The formal asymptotic expansion connects the homogenization problem with the ergodic cell problem. In this section, we use the formal asymptotic expansion to derive the effective Hamiltonian for multiscale homogenization problems.\\

We assume that the solution $u_{\e}(t,x)$ of (2) is developped as follows. (The evolutionary case (1) can be treated similarly, which we do not write here.) 
\begin{equation}
	u_{\e}(x)=\overline{u}(x)+ \sum^{k=d,n=N}_{k=1,n=1} \e^{n}_{k} u^{n}_{k} (x,\frac{x}{\e^1},...,\frac{x}{\e^N}), 
\end{equation}
where $u^n_k$ ($1\leq k\leq d$, $1\leq n\leq N$) are real valued functions defined in $(x,y^1,...,y^N)$$\in \Omega\times {\bf R^{d\times N}}$. \\
Then, the first derivatives are
\begin{equation}
	\frac{\p u_{\e}}{\p x_i}(x)=\frac{\p \overline{u}}{\p x_i}(x)+ \sum^{k=d,n=N}_{k=1,n=1} \sum^{N}_{m=1} \frac{\e^{n}_{k}}{\e^m_i} 
		\frac{\p u^n_k}{\p y^m_i}  (x,y^1,...,y^N)+ o(|\e|) \quad 1\leq i\leq d.
\end{equation}
If we put 
\begin{equation}
	w(x,y^1,...,y^N)=\sum^{k=d}_{k=1} \sum^{n=N}_{n=1} \frac{\e^{n}_{k}}{\e^1_1} 
		u^n_k(x,y^1,...,y^N)\quad 1\leq i\leq d, \quad 1\leq m\leq N,
\end{equation}
$$
	\gamma^m_i=\frac{\e^1_1}{\e^m_i} \quad 1\leq i\leq d,\quad 1\leq m\leq N, 
$$
and define $\Gamma^n$ ($1\leq n\leq N$) by (10), then 
\begin{equation}
	\frac{\p u_{\e}}{\p x_i}(x)=\frac{\p \overline{u}}{\p x_i}(x)+ \sum^{N}_{n=1} \Gamma^n \n_{y^n} w + o(|\e|)  \quad 1\leq i\leq d, 
\end{equation}
and by introducing (80) into (2) we get the ergodic problem (12)-(13). \\

\end{document}